\documentclass[12pt]{amsart}
\usepackage{amsmath,mathrsfs,amsthm}
\usepackage{txfonts,rotating,tikz}
\usepackage{esvect}

\usepackage[all]{xy}
\usepackage{nicefrac}

\usepackage[all]{xy}
\usepackage{color} 
\usepackage{hyperref}

\newcommand{\doi}[1]{\href{https://doi.org/#1}{\texttt{doi:#1}}}
\newcommand{\arxiv}[1]{\href{https://arxiv.org/abs/#1}{\texttt{arXiv:#1}}}

\newtheorem{theorem}{Theorem}[section]

\newtheorem{lemma}[theorem]{Lemma}
\newtheorem{proposition}[theorem]{Proposition}
\newtheorem{corollary}[theorem]{Corollary}

\newtheorem{conjecture}[theorem]{Conjecture}
\newtheorem{remark}[theorem]{Remark}
\newtheorem{assumption}[theorem]{Assumption}

\newtheorem{notation}[theorem]{Notation}
\newtheorem{definition}[subsection]{Definition}

\DeclareMathOperator{\Hom}{\mathrm{Hom}}

\DeclareMathOperator{\Res}{\mathrm{Res}}

\DeclareMathOperator{\Aut}{\mathrm{Aut}}
\DeclareMathOperator{\Out}{\mathrm{Out}}

\newcommand{\fh}{\mathfrak{h}}

\newcommand{\fb}{\mathfrak{b}}
\newcommand{\fg}{\mathfrak{g}}

\newcommand{\bbC}{{\mathbb {C}}}

\newcommand{\bbP}{{\mathbb {P}}}

\newcommand{\bbR}{{\mathbb {R}}}

\newcommand{\bbZ}{{\mathbb {Z}}}

\newcommand{\la}{\langle}
\newcommand{\ra}{\rangle}

\newcommand{\beq}{\begin{equation}}
\newcommand{\eeq}{\end{equation}}
\newcommand{\bt}{\begin{theorem}}
\newcommand{\et}{\end{theorem}}
\newcommand{\bde}{\begin{definition}}
\newcommand{\ede}{\end{definition}}
\newcommand{\bpr}{\begin{proposition}}
\newcommand{\epr}{\end{proposition}}
\newcommand{\ble}{\begin{lemma}}
\newcommand{\ele}{\end{lemma}}
\newcommand{\bco}{\begin{corollary}}
\newcommand{\eco}{\end{corollary}}
\newcommand{\bre}{\begin{remark}}
\newcommand{\ere}{\end{remark}}
\newcommand{\bpf}{\begin{proof}}
\newcommand{\epf}{\end{proof}}

\newcommand{\bc}{\mathbb{C}}

\newcommand{\bz}{\mathbb{Z}}

\begin{document}

\title{ Twisted conformal blocks and their dimension}

\author{Jiuzu Hong}
\address{Department of Mathematics, University of North Carolina at Chapel Hill, Chapel Hill, NC 27599-3250, U.S.A.}
\email{jiuzu@email.unc.edu}
\author{Shrawan Kumar}
\address{ 
Department of Mathematics, University of North Carolina at Chapel Hill, Chapel Hill, NC 27599-3250, U.S.A.}
\email{shrawan@email.unc.edu}

\maketitle
\begin{abstract}
Let $\Gamma$ be a finite group acting on a simple Lie algebra $\mathfrak{g}$ and acting on a $s$-pointed projective curve $(\Sigma, \vec{p}=\{p_1, \dots, p_s\})$ faithfully (for $s\geq 1$). Also, let an integrable highest weight module $\mathscr{H}_c(\lambda_i)$ of an appropriate twisted affine Lie algebra determined by the ramification at $p_i$  with a fixed central charge $c$ is attached to each $p_i$. 
We prove that the space of twisted  conformal blocks attached to this data
   is isomorphic to the space associated to a quotient group of $\Gamma$ acting on $\fg$ by diagram automorphisms and acting on a quotient of  $\Sigma$.   Under some mild conditions on ramification types, we prove that calculating the dimension of twisted conformal blocks can be reduced to the situation when $\Gamma$ acts on $\fg$ by  diagram automorphisms and covers of $\mathbb{P}^1$ with 3 marked points.  Assuming a twisted analogue of Teleman's vanishing theorem of Lie algebra homology, we  derive an analogue of the Kac-Walton formula and the Verlinde formula for general $\Gamma$-curves (with mild restrictions on ramification types). In particular, if the Lie algebra $\fg$ is not of type $D_4$, there are no restrictions on ramification types. 
\end{abstract}

\tableofcontents
\section{Introduction}
Wess-Zumino-Witten model is a type of two dimensional conformal field theory,
which associates to an algebraic curve with marked points and integrable highest weight modules of an affine Kac-Moody Lie algebra attached to the points, a finite dimensional vector space consisting of conformal blocks. The space of conformal blocks has many important properties including Propagation of Vacua and Factorization. It is also known that the sheaf of conformal blocks on the Deligne-Mumford stack of stable pointed curves is locally free.  The mathematical theory of WZW model was first established by Tsuchiya-Ueno-Yamada \cite{TUY} where all these properties were obtained. All the above properties are important ingredients in the proof of the celebrated Verlinde formula for the dimension of the space of conformal blocks (cf. \cite{Be, Fa, Kbook2, So1, V}).

One can replace algebraic curves by $\Gamma$-covers of curves for some finite group $\Gamma$, and let $\Gamma$ act on a simple Lie algebra $\fg$. Then, the theory of {\it twisted conformal blocks} can be similarly developed. It is related to the two dimensional orbifold conformal field theory in the literature \cite{BFH}, where Birke-Fuchs-Schweigert initiated this theory from the perspective of mathematical physics and conjectured an analogous  Verlinde formula for twisted conformal blocks in certain cases. In \cite{HK}, the authors obtained similar results as in \cite{TUY} for $\Gamma$-curves, including the properties of Propagation of Vacua and Factorization (under a technical assumption that $\Gamma$ stabilizes a Borel subalgebra of $\fg$; which is automatically satisfied if $\Gamma$ is cyclic), and we construct a flat projective  connection on the sheaf of twisted covacua on the Hurwitz stack of pointed smooth curves and we also prove local freeness of the sheaf of twisted covacua on the Hurwitz stack of stable pointed $\Gamma$-curves.  Earlier, similar results were obtained by Damiolini \cite{D1} under more restrictive conditions; in particular, where the marking points are unramified.

This paper is a continuation of our previous work \cite{HK}.  As our first main result of this work, in Theorem \ref{red_thm},  we prove that for any $\Gamma$-action on $\fg$, the dimension of the space of twisted conformal blocks is the same as the dimension of twisted conformal blocks attached to  $\tilde{\Gamma}$ acting on $\fg$ by diagram automorphisms and acting on a quotient curve $\tilde{\Sigma}$ of $\Sigma$, where $\tilde{\Gamma}$ is the quotient group of $\Gamma$ by the subgroup of elements acting on $\fg$ by inner automorphisms. In particular,  when $\Gamma$ acts on $\fg$ by inner automorphisms, the dimension of twisted conformal blocks is, in fact,   the same as the dimension of standard (nontwisted)  conformal blocks on the quotient curve, which can be computed by the usual Verlinde dimension formula,  cf.\, Corollary \ref{newcoro}.  Another application is given in Theorem \ref{thm_con_rank},  which asserts that if the quotient group $\tilde{\Gamma}$ is cyclic, then the sheaf of twisted conformal blocks on the Hurwitz stack of stable pointed $\Gamma$-curves is actually locally free of constant rank. Note that the sheaf of twisted conformal blocks
on the Hurwitz stack of stable pointed $\Gamma$-curves
 is proved to be locally free  in \cite[Theorem 8.9]{HK}. However, this stack may not be connected in general, and hence a priori it is unclear that the sheaf is of constant rank.

In Section 4,  we assume the group $\Gamma$ is cyclic. Under some restriction on ramification type at marked points,   in Theorem \ref{g_red_thm_2} we give a formula for the dimension of the twisted conformal blocks  in terms of the dimension  of the twisted conformal blocks 
for covers of $\mathbb{P}^1$ with $3$ marked points together with  the usual Verlinde numbers of higher genus.  This is achieved mainly by using the degeneration technique to create a node in $\Sigma$ and then using the Factorization Theorem thereby reducing the problem to a lower genus base curve $\bar{\Sigma}$ (cf. Lemma \ref{deg_lem1}). Further, by using a similar degeneration technique and the Factorization Theorem, we reduce the problem to a $\Gamma$-cover 
of $\mathbb{P}^!$ with only two ramified points 
(cf. Lemma \ref{dg_lem2}).

In Section 5,  we formulate a conjecture which is a  twisted analogue of Teleman's vanishing theorem for the Lie algebra homology (cf. Conjecture \ref{Teleman}). The proof of this conjecture 
 will appear in a separate work which is an on-going project by the authors \cite{HK2}. 
In Theorem \ref{KW_thm}, assuming the vanishing conjecture, we prove an analogue of the  Kac-Walton  formula for the dimension of twisted conformal blocks on covers of $\mathbb{P}^1$ by a cyclic group $\Gamma$ and  $\Gamma$ acting on $\fg$ by `standard' automorphisms (defined in Section \ref{sect_special_aut}). The main ingredient  in the proof of Theorem \ref{KW_thm} is the generalized Bernstein-Gelfand-Gelfand  resolution for twisted affine Kac-Moody Lie algebras, cf.\,Proposition \ref{BGG_comples}.

The first coauthor derived a Verlinde type formula for the trace of a diagram automorphism and defined twisted fusion rings in \cite{Ho1,Ho2}. These results (more specifically Theorem \ref{verlinde_Hong}) and Kac-Walton formula Theorem \ref{KW_thm} are two main ingredients in the proof of Theorem \ref{verlinde_genus_0}, which asserts that assuming the homology vanishing Conjecture \ref{Teleman}, there is a Verlinde type formula  for the dimension of twisted conformal blocks associated to covers of $\mathbb{P}^1$ with $3$ marked points  and standard automorphisms of $\fg$. 
Earlier, we expected a relation between the trace of diagram automorphism on a simple Lie algebra and the dimension of twisted conformal blocks for another related Lie algebra. Even though this explicit relationship is not exactly achieved, however the way we deduce the dimension of twisted conformal blocks associated to covers of $\mathbb{P}^1$ with $3$ marked points gives an indirect explanation of their relation. In particular, the formula in Theorem \ref{verlinde_Hong} and Theorem \ref{verlinde_genus_0} look fairly similar.

We finally combine all the above results and Conjecture  \ref{Teleman} to prove our second main 
 result of the paper: Theorem \ref{dimensionformula} determining the dimension of twisted conformal blocks in a fairly general setting (under some mild restriction on the ramification type only in the case of $\fg =D_4$). Specifically, Reduction Theorem  \ref{red_thm}; degeneration results Lemmas \ref{deg_lem1}  and \ref{dg_lem2} (resulting in Theorem  \ref{g_red_thm_2}); and  Theorem \ref{verlinde_genus_0}  for covers of $\mathbb{P}^1$ are the important ingredients in the proof of Theorem  \ref{dimensionformula}.

   Using the machinery of crossed modular categories, under the assumption that  $\Gamma$ stabilizes a Borel subalgebra of $\fg$ as in \cite{HK}, Deshpande-Mukhopadhyay \cite{DM}  deduced a Verlinde type formula for the dimension of twisted conformal blocks,  which is expressed in terms of S-matrices. The basic difference in their approach and ours is that we first of all reduce the problem to the standard automorphisms of $\fg$ and then we use the degeneration technique and the analogue of Kac-Walton formula to arrive at our dimension formula.

\vskip3ex
\noindent
{\bf Acknowledgements.} We thank Constantin Teleman for some helpful correspondences and conversations. The first author was partially supported by the NSF grant DMS-2001365 and the second
author was partially supported by the NSF grant DMS-1802328.
\section{Preliminaries}
\subsection{Kac-Moody theory}
\label{Kac_sect}

Let $\fg$ be a simple Lie algebra over $\mathbb{C}$.  Let $\sigma$ be an automorphism of order $m$ of $\fg$.  Let $\mathcal{K}$ be the field of Laurent series in the parameter $t$,  such that $\sigma(t)=\epsilon^{-1}t$ where $\epsilon=e^{ \frac{2\pi \mathrm{i}}{m} }$ and $\sigma$ acts on $\mathbb{C}$ trivially.      Let $\mathcal{O}$ be the field of formal power series in $t$. We now define a central extension  $\hat{L}(\fg,\sigma):=\fg(\mathcal{K})^\sigma\oplus \mathbb{C}C$ of $\fg(\mathcal{K})^\sigma$ under the bracket
  \begin{equation}  \label{eq1.1.1.4}
[x[P]+z C, x'[P'] +z' C] =
[x,x' ][PP'] +m^{-1}\Res_{t=0} \,\bigl(({dP})
P'\bigr) \langle x,x'\rangle C,  
 \end{equation}
for  $x[P],x'[P']\in \fg(\mathcal{K})^\sigma$, $z, z'\in\bc$; where
$\Res_{t=0}$ denotes the coefficient of $t^{-1}dt$ and $\langle\,,\,\rangle$ denotes the normalized invariant form on $\fg$
so that  the induced form on $\fg^*$ satisfies $\langle\theta,\theta\rangle =2$ for the highest root $\theta$ of $\fg$.

{\it Throughout the paper, we fix a positive integer (called the level) $c>0$.  We also fix an integer $s>0$ denoting the number of marked points.}

We use $D_{c,\sigma}$ to denote the set of highest weights of $\fg^\sigma$ which parametrizes the integrable highest weight modules of $\hat{L}(\fg, \sigma)$  of level $c$, where the level denotes the action of $C$,  see \cite[$\S$ 2]{HK}. When $\sigma$ is trivial, we also use $D_c$ to denote this set for brevity.    For each $\lambda\in D_{c,\sigma}$, we will denote by $(\mathscr{H}_c(\lambda), \rho_\lambda)$ (or for simplicity $\mathscr{H}_c(\lambda)$) the associated integrable highest weight module of $\hat{L}(\fg, \sigma)$ of level $c$.

There exists a  `compatible' Cartan subalgebra $\fh$ and a `compatible' Borel subalgebra $\fb \supset \fh$ of $\fg$  both stable under the action of $\sigma$ such that 
\begin{equation}  \label{eq1.1.1.0} \sigma=\tau \epsilon^{{\rm ad} h}, 
\end{equation}
where $\tau$ is a (possibly trivial) diagram automorphism of $\fg$ of order $r$ preserving $\fh$ and $\fb$, $\alpha (h)\in \bz$ for any root $\alpha$ of $\fg$ and $\epsilon^{{\rm ad} h}$ is the inner automorphism of $\fg$ such that for any root $\alpha$ of $\fg$, $\epsilon^{{\rm ad} h}$ acts on the root space $\fg_\alpha$ by the multiplication $\epsilon^{\alpha(h)}$, and $\epsilon^{{\rm ad} h}$ acts on $\fh$ by the identity. Here  $h$ is an element in $\fh^\tau$. In particular,  $\tau$ and $\epsilon^{{\rm ad}h}$ commute.  Moreover,  $r$ divides $m$, $\alpha(h)  \in \mathbb{Z}^{\geq 0}$  for  any positive  root  $\alpha$  of  $\fg^\tau$ and $\theta_0(h)\leq \frac{m}{r}$ where $\theta_0\in (\fh^\tau)^*$ denotes the following weight of $\fg^\tau$:
\[ \theta_0=\begin{cases}  
\text{ highest root of } \fg,   \text{ if } r=1\\
\text{ highest short root  of  } \fg^\tau, \text{ if } r>1 \text{ and }(\fg, r)\neq (A_{2n},2)\\
  2\cdot \text{highest short root} \text{ of } \fg^\tau, \,\text{ if } (\fg,r)=(A_{2n},2).  \end{cases}  \]

Let  $\hat{L}(\fg, \tau)$ denote the Lie algebra with the construction similar to  $\hat{L}(\fg, \sigma)$ where $\sigma$ is replaced by $\tau$, $m$ is replaced by $r$ and $\epsilon$ is replaced by $\epsilon^{\frac{m}{r}}$. There exists an isomorphism of Lie algebras (cf. \cite[Theorem 8.5]{Ka}):
  \begin{equation} \label{neweqn2.3.1}  \Psi_\sigma:   \hat{ L}(\fg, \tau)\simeq   \hat{L}(\fg, \sigma)  
\end{equation}
given by $C\mapsto C$ and
 $x[t^j]\mapsto    x[t^{\frac{m}{r}j+k  } ]$, for any $x$ an  $\epsilon^{\frac{m}{r}j}$-eigenvector of $\tau$, and $x$ also a $k$-eigenvector of ${\rm ad } \,h$.  Then, the isomorphism $\Psi_\sigma$ induces a bijection 
 \begin{equation}
 \label{weight_formula}
   D_{c,\sigma}\simeq D_{c,\tau}, \quad  \lambda \mapsto \bar{\lambda} . \end{equation}

\begin{remark}
{\rm The explicit description of $D_{c, \sigma}$  is given in \cite[Lemma 2.1]{HK} in terms of $\{n_{\lambda, i}\,|\,   i\in I(\fg^\tau)  \}$ defined there.  Also,  $\bar{\lambda}$ can be expressed in terms of numbers $a_i,a_i^\vee$ which can be read from \cite[p.54-55]{Ka} via \cite[Theorem 8.7]{Ka}. For the convenience of  readers, we would like to point out that there is  a typo in the formula for $\alpha_i^\vee$ in \cite[Theorem 8.7]{Ka}. The correct expression  is: $\alpha_i^\vee=1\otimes H_i + \frac{a_is_i r}{a_i^\vee m} K$. }
\end{remark}

\subsection{Twisted conformal blocks}
\label{sect_TCB}
Let $\phi\colon \Gamma\to {\rm Aut}(\fg)$ be a group homomorphism, where $\Gamma$ is a finite group and $ {\rm Aut}(\fg) $ is the group of Lie algebra automorphisms of $\fg$. Let $\Sigma$ be a reduced  projective $\Gamma$-curve over $\mathbb{C}$ such that no nontrivial element of $\Gamma$ fixes point-wise any irreducible component of $\Sigma$. 
{\it Unless otherwise stated, by a $\Gamma$-curve we will always mean such a $\Gamma$-curve.} For any $p\in \Sigma$, let $\Gamma_p$ be the stabilizer subgroup of $\Gamma$ at $p$. Then, $\Gamma_p$ is cyclic if $p$ is a smooth point of $\Sigma$. Let $\gamma_p$ be the ramification type at $p$, i.e., 
$\gamma_p$ is a generator of $\Gamma_p$ such that it acts on the tangent space $T_p\Sigma$ by the scalar multiplication $e^{\frac{2\pi \mathrm{i}}{|\Gamma_p|}}$.


Fix a tuple $\vec{p}=(p_1,p_2,\cdots, p_s)$ of distinct smooth points in $\Sigma$ such that any two distinct points are not in the same $\Gamma$-orbit.  Assume further that each irreducible component of $\bar\Sigma :=\Sigma/\Gamma$ contains at least one $\Gamma\cdot p_i$. 
{\it Then, such a $(\Sigma, \vec{p})$ is called a $s$-pointed $\Gamma$-curve.} For each $i$, let $\gamma_i$ be the ramification type at $p_i$.
Let $t_{i}$ be a $\gamma_{i}$-equivariant formal parameter at $p_i$, i.e., $\gamma_it_i=e^{-\frac{2 \pi i}{|\Gamma_{p_i}|}}t_i$. Let $\mathcal{K}_{p_i}$ denote the field $\mathbb{C}((t_i))$ of Laurent series and let $\hat{L}(\fg, \gamma_i)$ be the associated  twisted affine Lie algebra. In fact, it does not depend on the choice of equivariant $t_i$. We are also given a tuple $\vec{\lambda}=(\lambda_1,\lambda_2,\cdots, \lambda_s)$ of elements, where $\lambda_i\in D_{c, \gamma_i}$ for each $i$.  
Following \cite[$\S$ 3]{HK}, we define the following {\it space of covacua},
\begin{equation}
\label{t_cb_def}
  \mathscr{V}_{\Sigma, \Gamma, \phi}(\vec{p},  \vec{\lambda}):=   \frac{ \mathscr{H}_c(\lambda_{1})\otimes\cdots \otimes\mathscr{H}_c(\lambda_{s}) } { \fg[ \Sigma\backslash  \Gamma\cdot \vec{p} ]^\Gamma\cdot ( \mathscr{H}_c(\lambda_{1})\otimes\cdots \otimes\mathscr{H}_c(\lambda_{s})  ) }  ,  
\end{equation}
where $\fg[ \Sigma\backslash  \Gamma\cdot \vec{p} ]^\Gamma$ is the Lie algebra of $\Gamma$-equivariant maps from $\Sigma\backslash  \Gamma\cdot \vec{p}$ to $\fg$, and the action of $\Sigma\backslash  \Gamma\cdot \vec{p}$ on $\mathscr{H}_c(\lambda_{1})\otimes\cdots \otimes\mathscr{H}_c(\lambda_{k})$ is given by \cite[Definition 3.5]{HK}. It was proved in \cite{HK} that   twisted conformal blocks share similar properties with  usual conformal blocks, including Propagation, Factorization, WZW connection, etc. Some of these results are also proved in \cite{D1} under more restrictive assumptions.

\section{Reduction from general actions to diagram automorphisms}
\subsection{ A key lemma}
Let $G$ be a connected, simply-connected simple algebraic group over $\mathbb{C}$, and let $\Gamma$ be a finite group acting on $G$. Let $G_{\rm ad}$ denote the quotient of $G$ by its center. Then, $\Gamma$ acts on $G_{\rm ad}$ naturally. 
Let $\Sigma$ be a smooth projective connected curve over $\mathbb{C}$ with a faithful action of $\Gamma$. We regard $G_{\rm ad}$ as the group of inner automorphisms of $\fg$, which is a normal subgroup of the full automorphism group ${\rm Aut}(\fg)$. Hence,  ${\rm Aut}(\fg)$ acts on $G_{\rm ad}$ via conjugation.
Let ${\rm Out}(\fg)$ be the quotient group ${\rm Aut}(\fg)/G_{\rm ad}$.


\begin{lemma}
\label{key_lem}
Suppose that we are given two group homomorphisms $\phi, \psi\colon  \Gamma\to {\rm Aut}(\fg)$ such that $\phi\cdot \psi^{-1}:\Gamma \to  G_{\rm ad}$, $\gamma\mapsto \phi(\gamma)\psi(\gamma)^{-1}$.
For any $\Gamma$-stable affine open subset $\Sigma^*$ in $\Sigma$, if the action of $\Gamma$ on $\Sigma^*$ is free, then there exists a regular map $F\colon \Sigma^*\to G_{\rm ad}$ such that 
\[   F(\gamma\cdot p)= \phi(\gamma) F(p) \psi(\gamma)^{-1} ,   \quad  \quad   \forall p\in \Sigma^*, \gamma\in \Gamma.  \]
Note that $ \phi(\gamma)F(p)\psi(\gamma)^{-1} $ is well-defined as an element of  $G_{\rm ad}$, since $\phi(\gamma)\psi(\gamma)^{-1}\in G_{\rm ad}$ for any $\gamma\in \Gamma$.
\end{lemma}
\begin{proof}
Let $\mathscr{G}_\phi$ be the following (parahoric) Bruhat-Tits group scheme over $\bar{\Sigma}^*$, $\mathscr{G}_\phi:= \pi_*(\Sigma^*\times G_{\rm ad})^\Gamma$, where $\bar{\Sigma}^*=\Sigma^*/\Gamma$,  $\pi_*$ denote the Weil restriction from $\Sigma^*$ to  $\bar{\Sigma}^*$, and the upper subscript $^\Gamma$ denotes taking $\Gamma$-fixed point scheme under the action $\phi$ of $\Gamma$ via conjugation on $G_{\rm ad}$ ($\gamma\cdot g:= \phi(\gamma)g\phi(\gamma)^{-1}$ for $\gamma \in \Gamma$ and $g\in G_{\rm ad}$).

Recall that a $(\Gamma, G_{\rm ad}, \phi )$-bundle on $\Sigma^*$ is a right principal $\Gamma$-equivariant  $G_{\rm ad}$-bundle $\mathcal{P}$ on $\Sigma^*$ such that 
\[  \gamma( x\cdot g)=    \gamma(x) \cdot Ad_{\phi(\gamma ) }(g)   ,  \quad  \text { for any }  \gamma\in \Gamma,  g\in G_{\rm ad},  x\in \mathcal{P}.  \]   

We now construct a $(\Gamma, G_{\rm ad},\phi) $-bundle   structure $\mathscr{P}_\psi$
on $\Sigma^* \times G_{\rm ad}$ as follows ($\Gamma$ acts as a conjugation action on $G_{\rm ad}$ via $\psi$):
\begin{itemize}
\item $G_{\rm ad}$-bundle:  $(p,x) \cdot g=(p, x\cdot g)$, for any $p\in \Sigma^*$ and $g,x\in G_{\rm ad}$;
\item $\Gamma$-action:  $\gamma\cdot (p,x)= (\gamma\cdot p,  \psi(\gamma) x  \phi(\gamma)^{-1}   )$, for any $\gamma\in \Gamma$. 
\end{itemize}

The above bundle $\mathscr{P}_\psi$ by taking $\phi$ instead of $\psi$ is called the {\it trivial $(\Gamma, G_{\rm ad},\phi) $-bundle $\mathscr{P}^o$.}

  Since all the points in $\Sigma^*$ are unramified,   $\pi_*(\mathscr{P}_\psi )^\Gamma$ is a $\mathscr{G}_\phi$-bundle, cf.\,\cite[Proposition 2.9]{D2}.  
Now, we are in  position to apply Heinloth uniformization theorem (cf.\,\cite[Theorem 1]{He}), which asserts that as $\mathscr{G}_\phi$-bundles, $\pi_*(\mathscr{P}_\psi )^\Gamma$ is isomorphic to $\mathscr{G}_\phi$. Applying the inverse functor $\pi^*(\cdot)\times_{\pi^*( \mathscr{G}_\phi )  } (G_{\rm ad})_{\Sigma^*} $, we get an isomorphism of $(\Gamma, G_{\rm ad}, \phi)$-bundles  $\Phi\colon \mathscr{P}_\psi\simeq  \mathscr{P}^\circ$ (cf. \cite[Theorem 3.2]{D2}). From the construction of $\mathscr{P}_\psi$, we  see that there exists a regular map $F\colon \Sigma^*\to G_{\rm ad}$  such that 
\[ \Phi(p,x)= (p, F(p)x),  \quad \text{ for any $p\in \Sigma^*$ and $x\in G_{\rm ad}$} .  \]

By consideration of $\Gamma$-equivariance,  one can easily deduce that $F$ satisfies the desired property: 
\[  F(\gamma\cdot p)= \phi (\gamma) F(p) \psi(\gamma)^{-1}, \quad  \text{ for any $\gamma\in \Gamma$ and $p\in \Sigma^*$}. \] 

\end{proof}

\subsection{Reduction theorem} \label{reduction}

We consider the following setup: 

We are given  a group homomorphism $\phi:  \Gamma \to  {\rm Aut}(\fg) $ and  a projective irreducible smooth $s$-pointed $\Gamma$-curve $(\Sigma, \vec{p})$.   Let $\Sigma^*$ be the complement $\Sigma\backslash (\cup \Gamma\cdot p_i)$.
Let $\Gamma_0$ be the kernel of the map $P\circ \phi: \Gamma\to {\rm Out}(\fg)$, where $P: {\rm Aut}(\fg)\to {\rm Out}(\fg)$ is the projection map and ${\rm Out}(\fg)$ is the quotient group ${\rm Aut}(\fg)/{\rm Int}(\fg)$ (${\rm Int}(\fg)$ being the group of inner automorphisms of $\fg$).  Let $\tilde{\Gamma}$ be the quotient group $\Gamma/\Gamma_0$ and let $\tilde{\Sigma}$ be the quotient curve $\Sigma/ \Gamma_0$.  Let $\tilde{p}_i$ denote the image of $p_i$ in $\tilde{\Sigma}$, and let $\tilde{\Sigma}^*$ denote the complement $\tilde{\Sigma}\backslash (\tilde{\Gamma}\cdot \vec{\tilde{p} })$, where $\vec{\tilde{p}}= \{\tilde{p}_1, \dots, \tilde{p}_s\}$. Then, $(\tilde{\Sigma}, \vec{\tilde{ p}})$ is a $s$-pointed $\tilde{\Gamma}$-curve. Let $\tilde{\phi}_\iota $ be the composition of the following maps: 
\[  \tilde{\Gamma} \hookrightarrow   {\rm Out} (\fg) \xrightarrow{ \iota } {\rm Aut} (\fg) ,\]
 where $\iota$ is a group homomorphism such that the elements in $ {\rm Out} (\fg) $ act on $\fg$ by diagram automorphisms, which preserve a  pair $(\fb,\fh)$ and a pinning with respect to the pair $(\fb,\fh)$,
where $\fb$ is a fixed Borel subalgebra and $\fh$ is a Cartan subalgebra contained in $\fb$.  Let $\phi_\iota$ be the composition 
\[ \Gamma\to \tilde{\Gamma}\xrightarrow{\tilde{\phi}_\iota } {\rm Aut}(\fg).\]

For each $p_i$, choose a Borel subalgebra $\fb_i$ and a Cartan subalgebra $\fh_i$ contained in $\fb_i$ and both preserved by $\phi(\gamma_i)$ and satisfying the equation \eqref{eq1.1.1.0}, where $\gamma_i$ is the ramification type at $p_i$.
Let $\tau_i$ be the diagram automorphism part of $\phi(\gamma_i)$ with respect to this choice, i.e., the image of $\gamma_i$ under the analogue of $\phi_\iota$ with respect to the choice $(\fb_i, \fh_i)$.  

\begin{lemma}
\label{lem_kappa_innner}
For any $1\leq i\leq s$, there exists an inner automorphism $\kappa_i$ of $\fg$ such that $\kappa_i(\fb)=\fb_i$, $\kappa_i(\fh)=\fh_i$,  and
\begin{equation}
\label{iso_eq}
\tau_i= \kappa_i  \cdot  \phi_\iota(\gamma_i)  \cdot \kappa_i^{-1} .
\end{equation}
\end{lemma}
\begin{proof}
By an isomorphism theorem of semisimple Lie algebras, there exists an automorphism $\kappa'_i\in {\rm Aut}(\fg) $ such that $\kappa'_i(\fb)=\fb_i$, $\kappa'_i(\fh)=\fh_i$, and $\tau_i= \kappa'_i  \cdot  \phi_\iota(\gamma_i)  \cdot (\kappa'_i)^{-1}$. Let $D$ be the group of diagram automorphisms of $\fg$ preserving $\fb_i$ and $\fh_i$. It is well-known that $D \simeq {\rm Out}(\fg)$.  Thus, there exists an element $u\in D$ such that $\kappa_i :=u\kappa'_i$ is inner,  $\kappa_i(\fb)=\fb_i$ and $\kappa_i(\fh)=\fh_i$. Then, $\tau_i$ and  $ \kappa_i  \cdot  \phi_\iota(\gamma_i)  \cdot \kappa_i^{-1} $ are two elements in $D$. Note that $\tau_i$ and $\phi_\iota(\gamma_i)$ have the same image in ${\rm Out}(\fg)$. It follows that $\tau_i= \kappa_i  \cdot  \phi_\iota(\gamma_i)  \cdot \kappa_i^{-1} $.
\end{proof}

Given a tuple $\vec{\lambda}=(\lambda_1, \cdots, \lambda_s)$ of dominant weights with  $\lambda_i\in D_{c, \gamma_i}$, we get another tuple $\vec{\bar{\lambda}}=(\bar{\lambda}_1, \cdots, \bar{\lambda}_k )$ of dominant weights with $\bar{\lambda}_i\in D_{c, \tau_i}$ as described in (\ref{weight_formula}).  Via $\kappa_i$, we can  identify $D_{c,\tau_i}$ with  $D_{c, \tilde{\gamma}_i}$, where $\tilde{\gamma}_i = \phi_\iota(\gamma_i)$. We denote by $\tilde{\lambda}_i$ the element in $ D_{c, \tilde{\gamma}_i}$ corresponding to $\bar{\lambda}_i\in D_{c. \tau_i}$ under the identification $D_{c,\tau_i}\simeq D_{c, \tilde{\gamma}_i}$.

We  attach the space of twisted covacua $ \mathscr{V}_{\Sigma, \Gamma, \phi}(\vec{p},  \vec{\lambda})$ to  $(\Sigma, \vec{p})$ and $\phi\colon \Gamma\to {\rm Aut}(\fg)$. Similarly, we can also attach the space
$ \mathscr{V}_{\tilde{\Sigma}, \tilde{\Gamma}, \tilde{\phi}_\iota }(\vec{\tilde{p}},  \vec{\tilde{\lambda}})$ of twisted covacua  to the $s$-pointed $\tilde{\Gamma}$-curve $(\tilde{\Sigma}, \vec{\tilde{p}})$ and the group homomorphism $\tilde{\phi}_\iota\colon \tilde{\Gamma}\to {\rm Aut}(\fg)$. 
\begin{theorem}
\label{red_thm}
Assume that $\Sigma^* :=\Sigma\setminus(\Gamma\cdot \vec{p})$ does not contain any ramified points in $\Sigma$. Then, we have a natural isomorphism of vector spaces
\[  \mathscr{V}_{\Sigma, \Gamma, \phi} (\vec{p},  \vec{\lambda})\simeq  \mathscr{V}_{\tilde{\Sigma}, \tilde{\Gamma}, \tilde{\phi}_\iota }(\vec{\tilde{p}},  \vec{\tilde{\lambda}}) . \]
\end{theorem}
\begin{proof}
Let $\phi_\iota$ be the composition of the following maps:
\[  \Gamma \xrightarrow{\phi} {\rm Aut}(\fg) \xrightarrow{P}  {\rm Out} (\fg) \xrightarrow{ \iota } {\rm Aut} (\fg) , \]
where $\iota$ is as above.

By Lemma \ref{key_lem}, there exists $F\colon \Sigma^*\to G_{\rm ad}$ such that 
\begin{equation} 
\label{Thm_eq1}
 F(\gamma\cdot p)= \phi(\gamma) F(p) \phi_\iota(\gamma)^{-1} ,   \quad  \quad   \forall p\in \Sigma^*, \gamma\in \Gamma .\end{equation}
This gives rise to a Lie algebra homormophism $\Phi_F \colon \fg[\tilde{\Sigma}^* ]^{\tilde{\Gamma} } \to \fg[\Sigma^* ]^\Gamma $, given by 
\[  X\mapsto  {\rm Ad}_F (\pi^*X) , \quad  \text{for any } X\in \fg[\tilde{\Sigma}^* ]^{\tilde{\Gamma}},  \]
where $\pi^*X$ is the pull-back of the $\fg$-valued function $X$ on $\tilde{\Sigma}^*$, and ${\rm Ad}_F$ is the point-wise conjugation by $F$. One can check that $\Phi_F$ is an isomorphism. In fact, we construct its inverse map $\Psi_F$ as follows. For any $Y\in \fg[\Sigma^* ]^\Gamma$, it is easy to verify that ${\rm Ad}_{F^{-1}}(Y)\in \fg[\Sigma^* ]^{\Gamma, \phi_\iota }$, where $(\cdot)^{\Gamma, \phi_\iota }$ denotes the  $\Gamma$-invariants via the usual action of $\Gamma$ on $\Sigma^*$ and the action on $\fg$ via $\phi_\iota$. Then, ${\rm Ad}_{F^{-1}}(Y)$ descends to the desired element $\Psi_F(Y)\in \fg[\tilde{\Sigma}^* ]^{\tilde{\Gamma} }$.

Let $F_i$ be the image of $F$ in $G_{\rm ad}(\mathcal{K}_{p_i} )$.  Define  $\widehat{\rm Ad}_{F_{i}}:  \hat{L}(\fg, \tilde{\gamma}_i)\to \hat{L}(\fg, \gamma_i)$ as follows: 
\[  x[f]\mapsto  {\rm Ad}_{F_i} (x[f]) + \frac{1}{| \Gamma_{p_i}| }  {\rm Res}_{p_i} \langle F_i^{-1} d F_i , x[f]  \rangle C_i , 
\,\,\text{and}\,\, C_i\mapsto C_i,
 \]
where $ {\rm Ad}_{F_i}$ is the point-wise adjoint action,  $\langle,\rangle$ is the normalized invariant form on $\fg$, and $C_i$ is the canonical central element. Moreover,  $F_i^{-1} d F_i$ is the $\fg$-valued $1$-form,  which can be defined via an embedding $\rho\colon G_{\rm ad}\to GL(V)$. We regard $\mathcal{K}_{\tilde{p}_i }=\mathbb{C}((\tilde{t}_i))$ the subfield of $\mathcal{K}_{p_i}=\mathbb{C}((t_i))$, where $\tilde{t}_i=(t_i)^{e_i}$ and $e_i$ is the ramification index of $\pi\colon \Sigma\to \tilde{\Sigma}$ at $p_i$. Then, $ \widehat{\rm Ad}_{F_i} (x[f])\in \fg(\mathcal{K}_{\tilde{p}_i} )$.
 It is routine to check that $\widehat{\rm Ad}_{F_{i}}$ is a Lie algebra isomorphism. 

Let $ \oplus \hat{L}(\fg, \tilde{\gamma}_i) $ denote the direct sum of twisted affine Lie algebras  $\hat{L}(\fg, \tilde{\gamma}_i)$, and let $\overline{ \oplus }\hat{L}(\fg, \tilde{\gamma}_i)$ denote the quotient of $ \oplus \hat{L}(\fg, \tilde{\gamma}_i) $ by  the central  elements $C_i-C_j$ with $i\not=j$. Let $C$ denote the image of any $C_i$. Then, 
$\overline{ \oplus} \hat{L}(\fg, \tilde{\gamma}_i)$ has the  1-dimensional center $\mathbb{C} \cdot C$.  Similarly, we define the Lie algebra $\overline{ \oplus} \hat{L}(\fg,\gamma_i)$ with the canonical center $C$. Let $ \overline{\oplus}  \widehat{\rm Ad}_{F_{i}}\colon \overline{ \oplus} \hat{L}(\fg, \tilde{\gamma}_i)\to \overline{ \oplus} \hat{L}(\fg,\gamma_i)$ be the Lie algebra isomorphism induced from $\oplus  \widehat{\rm Ad}_{F_{i}}$.

We now consider the following diagram: 
\begin{equation}
\label{thm_diag}
\xymatrix{
\fg[\tilde{\Sigma}^* ]^{\tilde{\Gamma} } \ar[r]^{{\rm Loc}_{\vec{ \tilde{p} }}}  \ar[d]^{\Phi_F}  & \overline{\oplus} \hat{L}(\fg, \tilde{\gamma}_i) \ar[d]^{ \overline{\oplus}  \widehat{\rm Ad}_{F_{i}} } \\
\fg[\Sigma^* ]^\Gamma   \ar[r]^{{\rm Loc}_{\vec{p}}}  &   \overline{ \oplus} \hat{L}(\fg, \gamma_i),
}
\end{equation}
where ${\rm Loc}_{\vec{ \tilde{p} }}(Y) =\sum_i Y_{\tilde{p}_i }$ for any $Y\in \fg[\tilde{\Sigma}^* ]^{\tilde{\Gamma} } $, 
 and ${\rm Loc}_{\vec{p}}$ is defined similarly. By (\ref{Thm_eq1}), $F^{-1}dF$ is $(\Gamma,\phi_\iota)$-equivariant. It follows that the pairing $\langle F^{-1}dF, \pi^*Y \rangle$ is $\Gamma$-invariant 1-form on $\Sigma^*$. Hence, for any $q\in \Gamma\cdot p_i$,  the residue of $\langle F^{-1}dF, \pi^*Y \rangle$ at $q$ is equal to the residue at $p_i$.  Finally, the commutativity of the diagram (\ref{thm_diag}) follows from the following identity for any  $Y\in \fg[\tilde{\Sigma}^* ]^{\tilde{\Gamma} } $:  
 \[ \sum_i  \frac{1}{| \Gamma_{p_i}| }    {\rm Res}_{p_i}\langle F^{-1}dF, \pi^*Y \rangle  =  \frac{1} {|\Gamma|} \sum_{q\in \Gamma\cdot \vec{p}}    {\rm Res}_q\langle F^{-1}dF, \pi^*Y \rangle =0, \]
where the last equality follows from the Residue Theorem for $\langle F^{-1}dF, \pi^*Y \rangle$ on $\Sigma$
(cf. \cite[Chap. III, Theorem 7.14.2]{H}).

From the commutative diagram (\ref{thm_diag}), we have the following natural isomorphism:
\begin{equation}
\label{iso_1} 
 \mathscr{V}_{\Sigma, \Gamma, \phi}(\vec{p},  \vec{\lambda}) \simeq   \frac{   \mathscr{H}_c(\lambda_1) \otimes \cdots   \otimes   \mathscr{H}_c(\lambda_s)    }{   \fg[\tilde{\Sigma}^* ]^{\tilde{\Gamma} }   \cdot  (   \mathscr{H}_c(\lambda_1) \otimes \cdots   \otimes   \mathscr{H}_c(\lambda_s)    )  }  ,
 \end{equation}
where $ \mathscr{H}_c(\lambda_i)$ is regarded as a representation of $\hat{L}(\fg, \tilde{\gamma}_i)$ via the isomorphism $\widehat{\rm Ad}_{F_{i}}$.     

Recall the isomorphism $\Psi_{\gamma_i} \colon \hat{ L}(\fg, \tau_i)\simeq   \hat{L}(\fg, \gamma_i)$ from \eqref{neweqn2.3.1}.  It is an easy observation that $\Psi_{\gamma_i}=\widehat{\rm Ad}_{ F'_i} $, where $F'_i=t_i^{{\rm ad} h_i  } \in G_{\rm ad}( \mathcal{K}_{p_i} ) $,  $h_i\in \fh_i^{\tau_i}$ is determined by $\gamma_i$ as in Section \ref{Kac_sect}, and $t_i$ is the uniformizer in $\mathcal{K}_{p_i}$.  Moreover, by the equation \eqref{eq1.1.1.0}, 
$\phi(\gamma_i)= \epsilon_i^{{\rm ad}h_i  } \tau_i$, where $\epsilon_i$ is an $e_i$-th primitive root of unity.  Set $g_i=F'_i\cdot \kappa_i \cdot F_i^{-1} \in G_{\rm ad}(\mathcal{K}_{p_i})$, where $\kappa_i\in G_{\rm ad}$ is as in Lemma  \ref{lem_kappa_innner}
 thought of as an element of $G_{\rm ad}(\mathcal{K}_{p_i})$.  Then,  $g_i\in G_{\rm ad}(\mathcal{K}_{p_i} )^{\Gamma_{p_i}}$  since
\begin{align*}
 g_i(\gamma_i\cdot   ) &= F'_i(\gamma_i\cdot )\cdot \kappa_i \cdot F_i(\gamma_i\cdot )^{-1} \\
 &= t_i^{{\rm ad} h_i}(\cdot)  \epsilon_i^{{\rm ad} h_i } \kappa_i \phi_\iota(\gamma_i) F_i(\cdot)^{-1} \phi(\gamma_i)^{-1}  
 \\
 &= t_i^{{\rm ad} h_i}(\cdot)  \epsilon_i^{{\rm ad} h_i }  \tau_i   \kappa_i F_i(\cdot)^{-1} \phi(\gamma_i)^{-1}  \\
 &= t_i^{{\rm ad} h_i}(\cdot)    \phi(\gamma_i) \kappa_i F_i(\cdot)^{-1} \phi(\gamma_i)^{-1}  \\
 &=\phi(\gamma_i)  t_i^{{\rm ad} h_i}(\cdot)    \kappa_i F_i(\cdot)^{-1} \phi(\gamma_i)^{-1}   \\ 
&=\phi(\gamma_i)  g_i(\cdot)   \phi(\gamma_i)^{-1}  ,
    \end{align*}
where the second equality follows from (\ref{Thm_eq1}),  the third equality follows from (\ref{iso_eq}), and the fifth equality holds since $h_i\in \fh_i^{\gamma_i}$. 

 By a twisted anologue of Faltings' lemma  (cf.\,\cite[Proposition 10.2]{HK}), there exists an intertwining operator 
\[I_{g_i}\colon  (\mathscr{H}_c(\lambda_i), \rho_{\lambda_i} )\simeq ( \mathscr{H}_c(\lambda_i),   \rho_{\lambda_i}  \circ \widehat{\rm Ad}_{ g_i} )\]
 as isomorphisms of $  \hat{L}(\fg, \gamma_i) $-modules. This induces the following isomorphism of $  \hat{L}(\fg, \tau_i) $-modules
 \begin{align*}
 \tilde{I}_{g_i}\colon  (\mathscr{H}_c(\lambda_i),  \rho_{\lambda_i}   \circ  \widehat{\rm Ad}_{F_{i}}  )\simeq ( \mathscr{H}_c(\lambda_i),  \rho_{\lambda_i} \circ \widehat{\rm Ad}_{ g_i}  \circ    \widehat{\rm Ad}_{F_{i}} ) &\simeq (\mathscr{H}_c({\lambda}_i),  \rho_{\lambda_i}    \circ  \widehat{\rm Ad}_{F'_{i}} \circ \kappa_i  ) \\
&\simeq (\mathscr{H}_c(\bar{\lambda}_i),  \rho_{\bar{\lambda}_i  }\circ  \kappa_i  ),\,\,\text{by equation \eqref{weight_formula}}
\\
&\simeq  (\mathscr{H}_c(\tilde{\lambda}_i),  \rho_{\tilde{\lambda}_i  } ).
 \end{align*}
Therefore,  the operator (obtained from the above isomorphism $ \tilde{I}_{g_i}$ identifying  $  \hat{L}(\fg, \tilde{\gamma}_i) $ with  $  \hat{L}(\fg, \gamma_i) $ under $ \widehat{\rm Ad}_{F_{i}}$ as above):
\[ \otimes \bar{I}_{g_i}\colon   \mathscr{H}_c(\lambda_1) \otimes \cdots   \otimes   \mathscr{H}_c(\lambda_s) \simeq \mathscr{H}_c(\tilde{\lambda}_1) \otimes \cdots   \otimes   \mathscr{H}_c(\tilde{\lambda}_s) \] descends  to the following isomorphism
\begin{equation}
\label{iso_2}
\frac{   \mathscr{H}_c(\lambda_1) \otimes \cdots   \otimes   \mathscr{H}_c(\lambda_k)    }{   \fg[\tilde{\Sigma}^* ]^{\tilde{\Gamma} }   \cdot  (   \mathscr{H}_c(\lambda_1) \otimes \cdots   \otimes   \mathscr{H}_c(\lambda_s)    )  }  \simeq    \mathscr{V}_{\tilde{\Sigma}, \tilde{\Gamma}, \tilde{\phi}_\iota }(\vec{\tilde{p}},  \vec{\tilde{\lambda}}).
\end{equation} 
Combining the isomorphisms (\ref{iso_1}) and (\ref{iso_2}), we conclude the proof of this theorem.
\end{proof}

As a corollary of Theorem  \ref{red_thm}, we get the following result.

\begin{corollary} \label{newcoro} Let $\Gamma, \phi, \mathfrak{g}, \tilde{\Gamma}, (\Sigma, \vec{p})$ be as in the beginning of this Section \ref{reduction}. Assume that $\tilde{\Gamma} =(1)$. 
 We further assume that $\Gamma\cdot \vec{p}$ contains all the ramified points. Then, for any $\vec{\lambda}=(\lambda_1,\cdots, \lambda_{s})$ attached to $\vec{p}=(p_1,\cdots, p_{s})$ with  $\lambda_i\in D_{c,\gamma_i}$, 
  \[\dim \mathscr{V}_{\Sigma, \Gamma, \phi} (\vec{p}, \vec{\lambda}) = N_{\bar{g}} (\bar\lambda_1, \dots, \bar\lambda_s),\]
 where $\bar{g}$ is the genus of $\bar\Sigma :=\Sigma/\Gamma$, $\bar\lambda_i\in D_c$ is attached to $\lambda_i$ as in Section  \ref{reduction} and $N_{\bar{g}}(\bar\lambda_1, \dots, \bar\lambda_s)$ is the dimension of the untwisted conformal blocks attached to a genus $\bar{g}$ smooth irreducible curve and weights $ (\bar\lambda_1, \dots, \bar\lambda_s)$ attached to any distinct points. 
 
 For an explicit expression of  $N_{\bar{g}}(\bar\lambda_1, \dots, \bar\lambda_s)$ see \cite[Theorem 4.2.19]{Kbook2} or Theorem \ref{dimensionformula}. 
 
 In particular, the corollary applies for any non-simply laced $\mathfrak{g}$ (i.e., if $\fg$ is of type $B_\ell (\ell \geq 2), C_\ell (\ell\geq 2), F_4$ or $G_2$).
 \end{corollary}

 \subsection{An application}\label{application}
 
 We first recall the definition of stable $s$-pointed $\Gamma$-curves from \cite[Definition 8.1]{HK} (a variant of \cite[Definition 6.2.1]{BR}), where we call the same by `$\Gamma$-stable $s$-pointed $\Gamma$-curves'.
\begin{definition} \label{stablepointedcurves} {\rm A $s$-pointed $\Gamma$-curve $(\Sigma, \vec{p})$ (cf. $\S$\ref{sect_TCB})
is called {\it stable $s$-pointed $\Gamma$-curve} if $\Sigma$ is connected, $\bar{\Sigma} :=\Sigma/\Gamma$ is a  stable curve, i.e., it has at most  nodal singularity and the automorphism group of $(\bar{\Sigma}, \vec{\bar{p}})$ is finite (cf. \cite[Definition 2.1.1]{Kbook2}), where $\pi:\Sigma \to \bar{\Sigma}$ is the projection. Moreover, we require that for any node $q\in \Sigma$ and $\sigma \in \Sigma_q$, 
\begin{align*} \det (\dot{\sigma}) &= 1, \,\,\text{if $\sigma$ fixes the two branches at $q$}\\
&=-1 , \,\,\text{if $\sigma$ exchanges the two branches at $q$},
\end{align*}
where $\dot{\sigma}$ is the derivative of $\sigma$ acting on the Zariski tangent space $T_q(\Sigma)$.}
\end{definition}

We consider a stable $s$-pointed  $\Gamma$-curve $(\Sigma, \vec{p}=(p_1, \dots, p_s))$ of genus $g$  with marking data $\eta =\left((\Gamma_1, \chi_1),
 (\Gamma_2, \chi_2), \dots , (\Gamma_s, \chi_s)\right)$ (cf. \cite[Definition 8.7]{HK}).    By definition,  $\Gamma_i$ is the isotropy subgroup of $\Gamma$ at $p_i$.  We abbreviate $(\Gamma_i, \chi_i)$ by $\gamma_i$, where $\gamma_i$ is the generator of  $\Gamma_i$ such that its action on the tangent space $T_{p_i} (\Sigma)$ is via $e^{\frac{2 \pi \sqrt{-1}}{m_i}} Id$, where $m_i$ is the order of $\Gamma_i$.  Thus, the marking data $\eta$ can be identified with the ramification types $\vec{\gamma}=(\gamma_1,\gamma_2,\cdots, \gamma_s)$  at $\vec{p}$.  {\it We assume that $\Gamma\cdot \vec{p}$ contains all the ramified points in $\Sigma$.} 

\begin{remark}
{\rm Under the assumption that $\Gamma\cdot \vec{p}$ contains all the ramified points in $\Sigma$, at any nodal point $q\in \Sigma$, $q$ being unramified and stable, ${\rm det}(\dot \sigma)=1$,   $\sigma$ fixes the two branches for any $\gamma\in \Gamma_q$ and $\Gamma_q$ is cyclic (cf. \cite[Corollaire 4.3.3 and the comment after Definition 6.2.3]{BR}).  In this case, any stable $s$-pointed $\Gamma$-curve $(\Sigma, \vec{p})$ is exactly an admissible $s$-pointed $\Gamma$-cover in the sense of Jarvis-Kaufmann-Kimura \cite[Definition 2.1,2.2]{JKK}. The only difference is that, in our definition, stable $s$-pointed $\Gamma$-curves are connected, and admissible $s$-pointed $\Gamma$-covers defined in \cite{JKK} can be disconnected. }
\end{remark}
 
Let  $\overline{\mathscr{H}{M}}_{g, \Gamma,  \vec{\gamma}}$ be the Hurwitz stack of stable $s$-pointed $\Gamma$-curves of genus $g$ with marking data $\vec{\gamma}$, cf.\,\cite[$\S$ 8]{HK}. Then, $\overline{\mathscr{H}{M}}_{g, \Gamma,  \vec{\gamma}}$ is a proper and smooth Deligne-Mumford stack of finite type, cf.\,\cite[Theorem 8.8]{HK}.
We can attach the sheaf $\mathscr{V}_{g, \Gamma, \phi}(\vec{\gamma},  \vec{\lambda})$ of twisted covacua on $\overline{\mathscr{H}{M}}_{g, \Gamma,  \vec{\gamma}}$, where $\vec{\lambda}=(\lambda_1,\cdots,\lambda_s)$ with $\lambda_i\in D_{c, \gamma_i}$.  Then, $\mathscr{V}_{g, \Gamma, \phi}(\vec{\gamma},  \vec{\lambda})$ is locally free over $\overline{\mathscr{H}{M}}_{g, \Gamma,  \vec{\gamma}}$, cf.\,\cite[Theorem 8.9]{HK}. When $\Gamma$ is cyclic, $\overline{\mathscr{H}{M}}_{g, \Gamma,  \vec{\gamma}}$ is irreducible (\cite[Remark 8.11 (1)]{HK}). Thus, $\mathscr{V}_{g, \Gamma, \phi}(\vec{\gamma},  \vec{\lambda})$  is locally free of constant rank for cyclic $\Gamma$. 
For general $\Gamma$, $\overline{\mathscr{H}{M}}_{g, \Gamma,  \vec{\gamma}}$  could be disconnected. Nevertheless, we have the following theorem, which is an application of Theorem \ref{red_thm}.

 \begin{theorem}\label{thm_con_rank}
 With the notation as in Section \ref{reduction},   suppose that the quotient group $\tilde{\Gamma}$ of $\Gamma$ is cyclic. Then, the sheaf  $\mathscr{V}_{g, \Gamma, \phi}(\vec{\gamma},  \vec{\lambda})$ on   $\overline{\mathscr{H}{M}}_{g, \Gamma,  \vec{\gamma}}$ is locally free of constant rank. In particular, the theorem holds (without any assumption) for any $\fg$ of type other than $D_4$. 
 \end{theorem}
 \begin{proof} We freely follow the notation from Section \ref{reduction}. 
 Given any $s$-pointed smooth $\Gamma$-curve $(\Sigma, \vec{p})$ with ramification data $\vec{\gamma}=(\gamma_1,\cdots, \gamma_s)$ at $\vec{p}$, taking the quotient of $\Sigma$ by $\Gamma_0$ we get a smooth $s$-pointed $\tilde{\Gamma}$-curve with ramification data $\vec{\tilde{\gamma}}=(\tilde{\gamma}_1,\cdots, \tilde{\gamma}_s )$ at $\vec{\tilde{p}}$. Let $\tilde{g}$ be the genus of $\tilde{\Sigma}$.
 The Hurwitz stack of stable $s$-pointed $\tilde{\Gamma}$-curves with marking data $\vec{\tilde{\gamma}}$ is irreducible, since by assumption $\tilde{\Gamma}$ is cyclic.  By \cite[Theorem 8.9]{HK}, the sheaf $\mathscr{V}_{\tilde{g}, \tilde{\Gamma}, \tilde{\phi}_\iota}(\vec{\tilde{\gamma} },  \vec{\tilde{\lambda}})$ of twisted covacua on $\overline{\mathscr{H}{M}}_{\tilde{g}, \tilde{\Gamma},  \vec{\tilde{\gamma}} }$ is locally free of constant rank, where $\tilde{\phi}_\iota$ is the group action of $\tilde{\Gamma}$ on $\fg$ and $\vec{\tilde{\lambda}}$ is the $s$-tuple of dominant weights attached to $\vec{\tilde{p}}$ as in Theorem \ref{red_thm}.   By Theorem \ref{red_thm}, when $(\Sigma,\vec{p})$ is a smooth $s$-pointed $\Gamma$-curve, we have 
\[ \dim  \mathscr{V}_{\Sigma, \Gamma, \phi} (\vec{p},  \vec{\lambda})=\dim  \mathscr{V}_{\tilde{\Sigma}, \tilde{\Gamma}, \tilde{\phi}_\iota }(\vec{\tilde{p}},  \vec{\tilde{\lambda}}) . \]
This in particular implies that $\dim  \mathscr{V}_{\Sigma, \Gamma, \phi} (\vec{p},  \vec{\lambda})$ is constant along the smooth $s$-pointed $\Gamma$-curves $(\Sigma, \vec{p})$ in $\overline{\mathscr{H}{M}}_{g, \Gamma,  \vec{\gamma}}$.  By \cite[Theorem 8.9]{HK} again, the sheaf $\mathscr{V}_{g, \Gamma, \phi}(\vec{\gamma},  \vec{\lambda})$  is locally free.  To conclude the theorem, it suffices to show that every component of $\overline{\mathscr{H}{M}}_{g, \Gamma,  \vec{\gamma}}$ must contain a smooth $s$-pointed $\Gamma$-curve. Indeed this is true, as any stable $s$-pointed $\Gamma$-curve with nodal points admits a smoothing deformation (cf. \cite[Lemma 8.3 and  Proof of Theorem 8.9]{HK}). 

 \end{proof}

\section{Reduction via degenerations}

 In this section, we are in the same setup as in Section \ref{application} and  we further assume  that $\Gamma$ is cyclic of order $m$.

 Let $\bar{g}$ be the genus of $\bar{\Sigma}=\Sigma/\Gamma$. By the Riemann-Hurwitz formula when $\Sigma$ is a smooth $\Gamma$-curve, the genus $\bar{g}$ satisfies the following equation (cf. \cite[Corollary 2.4, Chap. IV]{H}):
\begin{equation}
\label{eqn3.0} 2g-2= |\Gamma| (2\bar{g}-2) +\sum_{i=1}^s \frac{|\Gamma|}{|\Gamma_i|}( |\Gamma_i| -1) .
\end{equation}


 
  \begin{lemma}
  \label{dim_const_lem} Let $(\Sigma, \vec{p})$ be a stable $s$-pointed smooth $\Gamma$-curve.
 Then, the  dimension of $ \mathscr{V}_{\Sigma, \Gamma, \phi}(\vec{p},  \vec{\lambda}) $ only depends on  $\phi, \bar{g}, \Gamma,  \vec{\gamma} =\{\gamma_1, \dots, \gamma_s\}, \vec{\lambda}$ and the level $c$. 
 \end{lemma}
 \begin{proof}
 By Riemann-Hurwitz formula (\ref{eqn3.0}), $g$ is determined by $\bar{g}, m$ and $\vec{\gamma}$. Thus, the lemma follows from Theorem \ref{thm_con_rank}. 
 \end{proof}
 Set (for fixed $\phi$ and $c>0$)
  \begin{equation}
  \label{verlinde_dim_eq}
  N_{\bar{g},\Gamma}(\vec{\gamma}, \vec{\lambda}  ) = \dim \mathscr{V}_{\Sigma, \Gamma, \phi}(\vec{p},  \vec{\lambda}) . \end{equation}

 \begin{lemma}
 \label{deg_lem1}
 Let $(\Sigma, \vec{p})$ be an irreducible   $s$-pointed smooth $\Gamma$-curve with ramification data $\vec{\gamma}$ such that $\Gamma\cdot \vec{p}$ contains all the ramified points in $\Sigma$. Assume that the quotient $\bar{\Sigma}$ has genus $\bar{g}\geq 1$  (in particular, $(\Sigma, \vec{p})$ is stable $\Gamma$-curve). Then, $(\Sigma, \vec{p})$ admits a degeneration to a stable $s$-pointed $\Gamma$-curve $(\Sigma', \vec{p'})$ (in particular,  $\Sigma'$ is connected) such that  the nodal points of $\Sigma'$ form a single $\Gamma$-orbit $\Gamma\cdot y$ and the action of $\Gamma$ on $\Gamma\cdot y$ is free. Moreover, 
$\Gamma\cdot \vec{p'}$ contains all the ramified points of $\Sigma'$.

If $\bar{g} \geq 2$ or if $\bar{g}=1$ and $\{\gamma_1, \dots, \gamma_s\}$ generate $\Gamma$, then  $\Sigma'$ can be taken to be irreducible. In any case, $\Sigma'/\Gamma$ is irreducible and hence we can take $ \vec{p'}$ to lie in an irreducible component of $\Sigma'$. 

 \end{lemma}
 \begin{proof}
 Let $\vec{\bar{p}}$ be the image of $\vec{p}$ in $\bar{\Sigma}$.  Then, the fundamental group of $\bar{\Sigma}\backslash \vec{\bar{p}}$ has the following presentation:
 \[  \big \{   \alpha_1,\beta_1,\cdots, \alpha_{\bar{g}}, \beta_{\bar{g}}, \eta_1,\cdots, \eta_s \,\mid\,   [\alpha_1,\beta_1]\cdots  [\alpha_{\bar{g}},  \beta_{\bar{g}}]\eta_1\cdots \eta_s    =1     \big  \}  , \]
 where $\eta_i $ represents the loop around the marked  point $\bar{p}_i$, and $\alpha_j,\beta_j$ represent loops arounds each handle of $\bar{\Sigma}$.
 The $\Gamma$-curve $\Sigma$ being irreducible gives rise  to a {\it surjective}  group homomorphism $f\colon  \pi_1( \bar{ \Sigma }  \backslash   \vec{\bar{p}})\to \Gamma$, where $\eta_i$ is mapped to $\gamma_i$ for each $1\leq i\leq s$.  In particular, we get $\gamma_1\gamma_2\cdots \gamma_s=1$ (since $\Gamma$ is cyclic by assumption; in particular, abelian).

 Let $(C', {\vec{\bar{p'}}})$ be a stable degeneration of $(\bar{\Sigma}, \vec{\bar{p}})$ with $C'$ irreducible and with one single node $\bar{x}$ (which is possible since $\bar{g}\geq 1$). Let $\tilde{C}'$ be the normalization of $C'$ with $\bar{x}^+,\bar{x}^-$ over $\bar{x}$.  Then, $\tilde{C}'$ is smooth and irreducible with genus $\bar{g}-1$. 
 Let $U$ be the complement 
 \[ C'\backslash \{\bar{x}, \bar{p}'_1,\cdots, \bar{p}'_s\}=\tilde{C}'\backslash \{ \bar{x}^+, \bar{x}^- , \bar{p}'_1,\cdots, \bar{p}'_s \}.\]  
 Then, the fundamental group of $U$ has the following presentation:
\[  \big \{   \alpha_1,\beta_1,\cdots, \alpha_{\bar{g}-1}, \beta_{\bar{g}-1}, \alpha^+, \alpha^- , \eta_1,\cdots, \eta_s \,\mid\,   [\alpha_1,\beta_1]\cdots  [\alpha_{\bar{g}-1},  \beta_{\bar{g}-1 }] \cdot \alpha^+ \alpha^- \eta_1\cdots \eta_s    =1     \big  \}  ,\]
where $\eta_i$ represents the loop around $\bar{p}'_i$, $\alpha^\pm$ represent the loops around $\bar{x}^{\pm}$, and $\alpha_j,\beta_j$ represent loops around each handle of $\tilde{C}'$.
We now construct a  group homomorphism $f'\colon  \pi_1(U)\to \Gamma$ such that $f'(\eta_i)=\gamma_i$ for any $1\leq i \leq s$, $f'(\alpha^+ )=f'(\alpha^-)=1$, and $f'(\alpha_j)=\gamma, f'(\beta_j)=\gamma^{-1}$ for each $1\leq j\leq \bar{g}-1$, where $\gamma$ is a generator of  the cyclic group $\Gamma$. Since $\gamma_1\cdots \gamma_s=1$, $f'$ is indeed a group homomorphism. 
The group homomorphism $f'$ gives rise to a $\Gamma$-bundle $\tilde{U} \to U$ with $\tilde{U}$ a smooth (but not necessarily connected) curve. By taking the unique smooth projective closure $\Sigma_{f'}\supset \tilde{U}$, we get  a smooth $s$-pointed $\Gamma$-cover $\pi\colon \Sigma_{f'} \to  \tilde{C}'$ with marked points $\vec{p'}$, such that the ramification data at ${\vec{p'}}$ is $\vec{\gamma}=(\gamma_1,\cdots, \gamma_s)$ and the ramification data above $\bar{x}^{\pm}$ is trivial. Let $y^\pm$ be a point above $\bar{x}^\pm$, chosen so that $y^+$ and $y^-$ are in the same component of the curve $\Sigma_{f'} $. Thus, $\pi^{-1}(\bar{x}^+)=\{\gamma^i\cdot y^+ \,|\, 0\leq i\leq m-1 \}$ and $\pi^{-1}(\bar{x}^-)=\{\gamma^i\cdot y^- \,|\, 0\leq  i\leq  m-1\, \}$ are free $\Gamma$-orbits.  By identifying $\gamma^i\cdot y^+ $ and $\gamma^{i+1}\cdot y^-$, for each $0\leq i\leq m-1$, we get a stable (in particular, connected)  $s$-pointed  $\Gamma$-curve $\Sigma'$ from $\Sigma_{f'}$ whose quotient by $\Gamma$ is exactly $C'$. Then, $(\Sigma', \vec{p'})$ is the desired stable $s$-pointed $\Gamma$-curve with nodal points $\pi^{-1}(\bar{x}^+)$. 
 \end{proof}

 \begin{lemma}
 \label{dg_lem2}
 Let $(\Sigma, \vec{p})$ be a stable $s$-pointed (irreducible) smooth $\Gamma$-cover of $(\mathbb{P}^1, \vec{\bar{p}})$ 
  (in particular, $s\geq 3$) such that $\Gamma\cdot \vec{p}$ contains all the ramified points in $\Sigma$ and has  ramification data $\vec{\gamma} = (\gamma_1, \dots, \gamma_s)$. 
 Suppose that $\gamma_1\gamma_2\cdots \gamma_t=1$ with both $ t, s-t\geq 2$. Then, the $\Gamma$-cover $\Sigma\to \bar{\Sigma}=\mathbb{P}^1$ degenerates to a stable $s$-pointed $\Gamma$-curve $(\Sigma', \vec{p'})$ whose quotient is a union of two projective lines  intersecting at a point $x$, such that above one projective line the ramification data is $(\gamma_1,\cdots, \gamma_t)$, and above another projective line the ramification data is $(\gamma_{t+1}, \cdots, \gamma_{s})$.  Moreover, the fiber over $x$ is a free $\Gamma$-orbit consisting of all the nodal points of $\Sigma'$. Further, 
$\Gamma\cdot \vec{p'}$ contains all the ramified points of $\Sigma'$. 

If $\{\gamma_1, \dots, \gamma_t\}$ generate $\Gamma$, then the curve over the first projective line can be taken to be irreducible. 
  \end{lemma}
 \bpf
 The fundamental group of $\mathbb{P}^1\backslash \{   \bar{p}_1,\bar{p}_2,\cdots, \bar{p}_s \}$ has a presentation:
 \[  \{  \eta_1,\eta_2,\cdots  ,\eta_s \,\mid \,  \eta_1\eta_2\cdots  \eta_s=1 \}    , \]
 where $\eta_i$ are loops around $\bar{p}_i$.
 
 The irreducible $\Gamma$-cover $\Sigma\to \mathbb{P}^1$ gives rise  to a surjective group homomorphism $f\colon \pi_1( \mathbb{P}^1 \backslash  \vec{\bar{p}})\to \Gamma$ such that $f(\eta_i)=\gamma_i$.  Let $U_1$ be $\mathbb{P}^1\backslash \{\bar{p}_1,\cdots, \bar{p}_t\} $ and let $U_2$ be  $\mathbb{P}^1\backslash \{\bar{p}_{t+1},\cdots, \bar{p}_s\}$. 
 For each $k=1,2$, we construct a group homomorphism 
 $f_k\colon \pi_1(U_k) \to  \Gamma$ such that 
\[   f_1(\eta_i)=\gamma_i, \text{ for any }1\leq i\leq t, \quad   f_2(\eta_j)= \gamma_j, \text{ for any } t+1\leq j\leq s.   \]
Observe that $f_1$ and $f_2$ are group homomorphisms, since by assumption $\gamma_1\gamma_2\cdots \gamma_t=1$. For each $k=1,2$, let $( \Sigma_{f_k}, \vec{p}_k)$ be the unique  smooth $\Gamma$-cover of $\mathbb{P}^1$ associated to $f_k$ ($\Sigma_{f_k}$ could be disconnected), such that $\vec{p}_k$ has the  ramification data $(\gamma_1,\cdots, \gamma_t)$ and  $(\gamma_{t+1}, \cdots, \gamma_{s})$ for $k=1,2$ respectively.   
Fix any (unramified) $\Gamma$-orbits  $\Gamma\cdot x_1\subset \Sigma_{f_1}$ and $\Gamma\cdot x_2\in \Sigma_{f_2}$
over points in $U_1$ and $U_2$ respectively.  We glue $\Sigma_{f_1}$ and $\Sigma_{f_2}$ along  any $\Gamma$-equivariant map between $\Gamma\cdot x_1$ and $\Gamma\cdot x_2$. Since $f$ is surjective, we get a connected $s$-pointed $\Gamma$-curve $(\Sigma', \vec{p}_1,\vec{p}_2)$ whose quotient is a union of two projective lines intersecting at a point $x$, with marked points $(\bar{p}_1,\cdots, \bar{p}_t, \bar{p}_{t+1},\cdots, \bar{p}_s)$.  This  $s$-pointed $\Gamma$-curve has the desired properties.

 \epf
 
 
 Let $\Gamma, \phi, \Sigma, \vec{p}=(p_1,\cdots, p_s), \vec{\gamma}= (\gamma_1, \dots, \gamma_s)$ and   $\vec{\lambda}=(\lambda_1,\cdots, \lambda_s)$ be as in Section  \ref{sect_TCB}. Assume further that $\Sigma$ is smooth and it  has an irreducible component $\Sigma^o$ such that each $p_i$ belongs to $\Sigma^o$. Moreover, $\bar{\Sigma} = \Sigma/\Gamma$ is irreducible.  Let $\Gamma^o$ be the subgroup of $\Gamma$ stabilizing $\Sigma^o$. Then, $\Gamma^o, \phi^o:=\phi_{|\Gamma^o}, \Sigma^o, \vec{p}, \vec{\gamma}$ and   $\vec{\lambda}$ also satisfy the assumptions of  Section  \ref{sect_TCB} (observe that $\Gamma_{p_i}= \Gamma^o_{p_i}$). 
 
 We have the following reduction lemma.
\begin{lemma}
\label{red_lem}
With the assumption as above, we have: 
\begin{enumerate}
  \item  There exists an isomorphism  $\beta:\Gamma\times_{\Gamma^o} \Sigma^o\simeq \Sigma$ of $\Gamma$-curves given by $[\gamma, x] \mapsto\gamma\cdot x$. 
\item  There exists an isomorphism of vector spaces  $ \mathscr{V}_{\Sigma, \Gamma, \phi}(\vec{p},  \vec{\lambda}) \simeq  \mathscr{V}_{\Sigma^o, \Gamma^o, \phi^o }(\vec{p},  \vec{\lambda}).  $  In other words, 
\[\dim  \mathscr{V}_{\Sigma, \Gamma, \phi}(\vec{p},  \vec{\lambda})=\dim \mathscr{V}_{\Sigma^o, \Gamma^o, \phi^o}
(\vec{p}, \vec{\lambda}) .  \]
\end{enumerate}
\end{lemma}
 \begin{proof}
 For part (1), clearly the map $\beta$  is surjective (since $\bar\Sigma$ is irreducible) and $\Gamma$-equivariant.  The injectivity follows from the definition of $\Gamma^o$ since, for any $\gamma\in \Gamma$, $\gamma\Sigma^o\cap \Sigma^o \not=\emptyset$ if and only if $\gamma\in \Gamma^o$ (this uses the smoothness of $\Sigma$). 
 
For part (2), it suffices to check that the restriction map 
 \[ {\rm Res}\colon \fg[ \Sigma\backslash \Gamma\cdot \vec{p} ]^\Gamma\to \fg[\Sigma^o \backslash \Gamma^o \cdot \vec{p} ]^{\Gamma^o}\]
 is an isomorphism:  For any $\Gamma$-equivariant map $X\colon \Sigma\backslash \Gamma\cdot \vec{p} \to \fg$, if $X$ vanishes on the component $\Sigma^o\backslash \Gamma^o \cdot \vec{p}$, then, by $\Gamma$-equivariance,  $X$ vanishes everywhere. Thus, the restriction map  ${\rm Res}$ is injective. 
  
  For any $\Gamma^o$-equivariant map $Y\colon  \Sigma^o \backslash \Gamma^o \cdot \vec{p}\to \fg$,  construct an extension $\tilde{Y}\colon \Sigma\backslash \Gamma\cdot \vec{p}\to \fg$ given by $\tilde{Y}(q)= \phi_{\gamma^{-1}} (  Y(\gamma\cdot q) )$ for any $\gamma$ such that $\gamma\cdot q\in \Sigma^o\backslash\Gamma^o\cdot \vec{p}$,  where $\phi_{\gamma^{-1}}$ is the automorphism of $\fg$ associated to $\gamma^{-1}$. One can check easily that $\tilde{Y}$ is a well-defined $\Gamma$-equivariant regular map. Thus, the restriction map ${\rm Res}$ is an isomorphism.  
  
 \end{proof}



\begin{notation}
\label{verlinde_notation}
\begin{enumerate}
{\rm \item Let $N_{g}(\vec{\lambda})$ denote the dimension of the space of (untwisted) conformal blocks attached to an irreducible smooth projective curve $C$ of genus $g$ and $\vec{\lambda}=(\lambda_1,\cdots, \lambda_s)$ at $s$-points in $C$ with $\lambda_i\in D_c$.
\item  For an integer $m\geq 1$, let $\Gamma_m=\langle \gamma \rangle$ (cyclic group of order $m$) act on $\mathbb{P}^1$ by $\gamma \cdot z=e^{ \frac{2\pi \mathrm{i}}{m} } z $ for  $z\in \mathbb{P}^1$ and $\phi:\Gamma_m\to \Aut \fg$. Let $(\lambda, \mu, \nu)$ be a set of dominant weights, such that $\lambda\in D_{c, \gamma}, \mu\in D_{c, \gamma^{-1}}$ and $\nu\in D_c$ attached to the points $(0, \infty, 1)$ respectively. We denote by $N_{\phi}(\gamma;\lambda,\mu,\nu)$ the dimension of the twisted conformal blocks attached to this data. }

\end{enumerate}
\end{notation}
  
It is well-known that $N_{g}(\vec{\lambda})$ can be computed by the usual Verlinde formula  (cf. \cite [Theorem 4.2.19]{Kbook2}).   By the reduction theorem in Theorem \ref{red_thm}, the computation of  $N_\phi(\gamma;\lambda,\mu,\nu)$ can be reduced to the case when $\gamma$ acts via a diagram automorphism of $\fg$. In fact, by the same reason it suffices to assume that $\gamma$ acts on $\fg$ via a  standard automorphism in the sense of Section \ref{sect_special_aut}.
In Section \ref{twisted_Verlinde_sect}, we will prove a Verlinde type  formula  for $N_\phi(\gamma;\lambda,\mu,\nu)$ when $\gamma$ acts on $\fg$ via a  standard automorphism. 
  \begin{lemma}
  \label{cover_lem}
Let $\Gamma$ be a cyclic group of order $m\geq 2$. Then, any irreducible smooth $\Gamma$-cover $\Sigma$ of $\mathbb{P}^1$ with two branched points in $\mathbb{P}^1$ is isomorphic to $\pi: \mathbb{P}^1\to \mathbb{P}^1$ given by $z\mapsto z^m$. 
  \end{lemma}
  \begin{proof}
 Let $g$ be the genus of $\Sigma$. Let $\bar{p}_1,\bar{p}_2$ be the branched points in $\mathbb{P}^1$ with ramification indices $e_1$ and $e_2$.   By Riemann-Hurwitz formula (cf. the identity \eqref{eqn3.0}),  
 \[  2g-2= -2m+ \frac{m }{e_1}( e_1-1 )+  \frac{m}{e_2}(e_2-1) = -\frac{m}{e_1}-\frac{m}{e_2} \leq -2.\]
 Thus, $g=0$, and $e_1=e_2=m$. It is easy to see that such a $\Gamma$-cover $\Sigma$ over $\mathbb{P}^1$ is isomorphic to $\pi: \mathbb{P}^1\to \mathbb{P}^1$ given by $z\mapsto z^m$.
 \end{proof}
 

 The following theorem reduces the problem of calculating the dimension of twisted conformal blocks to that of the classical Verlinde numbers together with $N_\phi$ as in Notation \ref{verlinde_notation}  
  \bt
 \label{g_red_thm_2}
 Let $(\Sigma, \vec{p})$ be an irreducible smooth $s$-pointed $\Gamma$-curve,  where $s\geq 1$ if $\bar{g}\geq 1$ and $s\geq 3$ if $\bar{g}=0$ (so that $(\Sigma, \vec{p})$ is a stable $s$-pointed $\Gamma$-curve) and $\Gamma$ is any finite cyclic group. 
 Assume that $\Gamma\cdot \vec{p}$ contains all the ramified points, and that we can write $\vec{p}=(p_1,\cdots, p_s)$ so that  $p_1,\cdots p_{2a}$ are ramified ($a \geq 0$) and $p_{2a+1},\cdots, p_{s}$ are unramified. Assume that $\gamma_{2k-1}\gamma_{2k}=1$ for each $1\leq k\leq a$. 
Let $\vec{\lambda}=(\lambda_1,\cdots, \lambda_{2a})$ be attached to $(p_1,\cdots, p_{2a})$ with each $\lambda_i\in D_{c,\gamma_i}$, and $\vec{\mu}=(\mu_1,\cdots, \mu_b)$ attached to $(p_{2a+1}, \cdots,p_{s})$ with each $\mu_j\in D_c$, where $s=2a+b$. 
 Then, we have the following formula:
 \[  N_{\bar{g},\Gamma}(\vec{\gamma};\vec{\lambda}, \vec{\mu})= \dim \mathscr{V}_{\Sigma, \Gamma, \phi} (\vec{p}, \vec{\lambda}, \vec{\mu}) = 
 \sum_{\vec{\nu}} \left (  \prod_{k=1}^a N_\phi(\gamma_{2k-1}; \lambda_{2k-1},\lambda_{2k}, \nu_k) \right ) \cdot N_{\bar{g}}(\vec{\mu},\vec{\nu}^*  )  , \]
 where the summation is over $\vec{\nu}=(\nu_1,\cdots, \nu_a)$ with $\nu_i \in D_c$. Here, $N_\phi(\gamma_{2k-1}; \lambda_{2k-1},\lambda_{2k}, \nu_k) $ and $N_{\bar{g}}(\vec{\mu},\vec{\nu}^*  )$ are defined in Notation \ref{verlinde_notation}. 
  \et
 \bpf
 We prove the theorem by reducing the problem for $\bar{g}$ to that of $\bar{g}-1$. So, assume that $\bar{g} \geq 1 $.    
 
 By Lemma \ref{deg_lem1},  there exists a stable $s$-pointed $\Gamma$-curve $(\Sigma', \vec{p'})$ in $\overline{\mathscr{H}M}_{g, {\Gamma},  \vec{\gamma}}$ (in particular,  $\Sigma'$ is connected)  and $\Sigma'$ has a $\Gamma$-free single $\Gamma$-orbit of  nodal points.  Moreover,  $\Gamma \cdot \vec{p'}$ contains all the ramified points of $\Sigma'$.  
 Further,  by part (1) of Lemma \ref{red_lem},   the normalization of $\Sigma'$ at nodal points is isomorphic to $\Gamma \times_{\Gamma^1} \Sigma_{1}$, where $\Sigma_{1}$ is an irreducible smooth projective  $\Gamma^{1}$-curve for some subgroup $\Gamma^{1}$ of $\Gamma$.   Applying the Factorization Theorem \cite[Theorem 5.4]{HK} and part (2) of Lemma \ref{red_lem}, we get:
 \begin{equation}
 \label{eqn3.2}
 N_{\bar{g}, \Gamma } (\vec{\gamma}; \vec{\lambda}, \vec{\mu}) = \sum_{\eta_1\in D_c} \,N_{\bar{g}-1, \Gamma^1}(\vec{\gamma}, 1, 1; \vec{\lambda}, \vec{\mu}, \eta_1, \eta_1^*).
 \end{equation} 
Thus, inducting on $\bar{g}$ and keep using the Factorization Theorem, Lemma \ref{deg_lem1} and Lemma \ref{red_lem}, we get
  \begin{equation}
 \label{eqn3.3}
  N_{\bar{g}, \Gamma} (\vec{\gamma}; \vec{\lambda}, \vec{\mu}) = \sum_{\vec{\eta}\in (D_c)^{\bar{g} }} \,N_{0, \Gamma'}(\vec{\gamma}, \vec{1}_{2\bar{g}}; \vec{\lambda}, \vec{\mu},  \vec{\eta}, \vec{\eta}^* ),
 \end{equation}  
for some subgroup $\Gamma'$ of $\Gamma$, where $\vec{\eta}=(\eta_1,\cdots, \eta_{\bar{g}})\in (D_c)^{\bar{g}}$ and $\vec{\eta}^* :=(\eta_1^*,\cdots, \eta^*_{\bar{g}}) \in   (D_c)^{\bar{g}}$. We emphasize here that $N_{0, \Gamma'}(\cdot)$ denotes the dimension of the space of twisted conformal blocks attached to an irreducible smooth projective $\Gamma'$-cover of $\mathbb{P}^1$. 

Similarly, keep using the Factorization Theorem, Lemma \ref{dg_lem2} and Lemma \ref{red_lem} for the pair $(p_{2k-1}, p_{2k})$ of points with $1\leq k\leq a$, we get (using the equation \eqref{eqn3.3}: 
\begin{equation}
\label{eq_3.9_1}
  N_{\bar{g}, \Gamma}(\vec{\gamma};\vec{\lambda}, \vec{\mu})=\sum_{\vec{\nu}, \vec{\eta}} \left ( \prod_{k=1}^a N_{0, \Gamma_{k}}(\gamma_{2k-1}, \gamma_{2k},1; \lambda_{2k-1},\lambda_{2k}, \nu_k) \right ) \cdot N_{0, \Gamma'' }(\vec{1}_{2\bar{g}+b+a}; \vec{\mu}, \vec{\eta}, \vec{\eta}^*, \vec{\nu}^*  )  , \end{equation}
for some subgroups $\Gamma''$ and $\Gamma_k$ of $\Gamma$ for each $1\leq k\leq a$,  where the summation is over $\vec{\nu}=(\nu_1,\cdots, \nu_a)$ and $\vec{\eta}=(\eta_1,\cdots, \eta_{\bar{g}}) $ with $\nu_i, \eta_j \in D_c$. 

Note that any \'etale $\Gamma''$-cover $\Sigma''$
over $\mathbb{P}^1$ is isomorphic to $\mathbb{P}^1\times \Gamma''$. Since $\Sigma''$ is irreducible, it
 follows that $\Gamma''=1$.  Then, by  Notation \ref{verlinde_notation}, we have 
\[  N_{0, \Gamma''}(\vec{1};\vec{\mu}, \vec{\eta}, \vec{\eta}^*, \vec{\nu}^* )  = N_0(\vec{\mu}, \vec{\eta}, \vec{\eta}^*, \vec{\nu}^*  )  . \]

By Lemma \ref{cover_lem}, any irreducible $\Gamma_k$-cover over $\mathbb{P}^1$ with ramification data $(\gamma_{2k-1}, \gamma_{2k})$  is isomorphic to a standard $\Gamma_{k}$-action on $\mathbb{P}^1$ as in Lemma \ref{cover_lem}.
 Then, by  Notation \ref{verlinde_notation},  for each $k$ we have 
\[  N_{0, \Gamma_k}(\gamma_{2k-1}, \gamma_{2k},1; \lambda_{2k-1},\lambda_{2k}, \nu_k) =N_\phi(\gamma_{2k-1}; \lambda_{2k-1},\lambda_{2k}, \nu_k)  .  \]
Thus, from (\ref{eq_3.9_1}) we get 
 \begin{align}
  N_{\bar{g}, \Gamma}(\vec{\gamma};\vec{\lambda}, \vec{\mu})&=\sum_{\vec{\nu}, \vec{\eta}} \left ( \prod_{k=1}^a N_\phi(\gamma_{2k-1}; \lambda_{2k-1},\lambda_{2k}, \nu_k) \right ) \cdot N_0(\vec{\mu}, \vec{\eta}, \vec{\eta}^*, \vec{\nu}^*  ) \\
  &=\sum_{\vec{\nu}} \left ( \prod_{k=1}^a N_\phi(\gamma_{2k-1}; \lambda_{2k-1},\lambda_{2k}, \nu_k) \right ) \cdot \sum_{ \vec{\eta}} N_0(\vec{\mu}, \vec{\eta}, \vec{\eta}^*, \vec{\nu}^*  )\\
  &=\sum_{\vec{\nu}} \left ( \prod_{k=1}^a N_\phi(\gamma_{2k-1}; \lambda_{2k-1},\lambda_{2k}, \nu_k) \right ) \cdot N_{\bar{g}}(\vec{\mu}, \vec{\nu}^*  ),   \end{align}
  where the lasts equality follows from factorization of fusion rules for conformal blocks in untwisted setting, cf.\,\cite[Corollary 3.5.10 (a)]{Kbook2}.  This concludes the proof of the theorem. 
 \epf
 
\begin{remark} {\rm Assume that $\Gamma\simeq \mathbb{Z}/2\mathbb{Z}$.  Let $(\Sigma, \vec{p})$ be a stable smooth $s$-pointed $\Gamma$-curve such that $\Gamma\cdot \vec{p}$ contains all the ramified points in $\Sigma$. By Riemann-Hurwitz formula (\ref{eqn3.0}), there are even number of ramified points in $\Sigma$. Thus, up to ordering we can always write 
 \[ \vec{p}=(p_1,\cdots, p_{2a}, p_{2a+1},\cdots, p_{s}), \]
  so that $p_1,\cdots, p_{2a}$ are ramified and $(p_{2a+1}, \cdots, p_{a})$ are unramified.  Then, by Theorem \ref{g_red_thm_2},  the dimension of  the space of twisted conformal blocks attached to any ramification data can be reduced to compute $N_\phi(\gamma; \lambda, \mu,\nu)$ and $N_{\bar{g}}(\vec{\lambda})$. }
\end{remark}
\begin{remark}
{\rm Assume that $\Gamma\simeq \mathbb{Z}/3\mathbb{Z}$. We have an elliptic curve $E$ over $\mathbb{P}^1$ as a $\Gamma$-cover. The ramification type of $E$ is $(\gamma, \gamma, \gamma)$, where $\gamma$ is a generator of $\Gamma$.  In this case, Theorem \ref{g_red_thm_2}
 is not applicable since $\gamma^2\neq 1$. }
\end{remark}

  \section{Kac-Walton formula for twisted conformal blocks}

\subsection{Standard automorphisms}
\label{sect_special_aut}
An automorphism $\sigma$ of $\fg$ is called $\bold{special}$ if $\sigma$ is a  diagram automorphism
 (which includes the identity automorphism), or an order 4 automorphism of $\fg$ when $\fg$ is of type $A_{2n}$, which is defined as follows. Let $e_i,f_i,h_i, i=1,\cdots, 2n$,  be the set of Chevalley generators. The automorphism $\sigma$ of $\fg$ is defined such that 
\begin{equation}
\label{aut_def1}
\begin{cases}
\sigma(e_i)=e_{\tau(i)},   \quad  \text{ if } i\not= n, n+1;    \\
 \sigma(e_i)=\sqrt{-1} e_{\tau(i)},   \quad \text{ if } i\in \{n, n+1\} ; \\
  \sigma(f_\theta)=f_\theta ,
\end{cases},
\end{equation}
where $\theta$ is the highest root of $\fg$ and $\tau$ is the nontrivial diagram automorphism.  In fact, we can write 
\beq  
\label{as_auto}
 \sigma= \tau \sqrt{-1}^{{\rm ad }h },\eeq
 where $h\in \fh$ is such that $\alpha_i(h)=0$ if $i\not=n, n+1$ and $\alpha_i(h)=1$ if $i=n,n+1$. 

We call $\sigma$ to be a {\it standard special} automorphism (or simply a {\it standard} automorphism) if $\sigma$ is the identity automorphism or 
 a nontrivial diagram automorphism when $\fg$ is not of type $A_{2n}$ or $\sigma$ is the order 4 special automorphism as above. 
(Observe that for  a standard automorphism,   $\fg^\sigma$ is the same as $\mathring{\fg}$ as  defined in \cite{Ka}.) So, the only difference between special and standard automorphism is that we exclude the nontrivial diagram automorphism of $A_{2n}$. 

The following table describes the fixed point Lie algebra for all the nontrivial special automorphisms, cf.\,\cite[$\S$ 2.1]{BH}:
 \begin{equation}
\label{Fix_table}
 \begin{tabular}{|c  | c | c |c |c |c | c| c| c|c |c|c|c|c|c |c ||} 
 \hline
$(\fg, m)$   & $(A_{2n-1}, 2 ) $  &  $(A_{2n}, 4)$  &   $(A_{2n}, 2)$  &  $(D_{n+1} , 2  ) $  &           $ (D_4,   3)$  &  $ (E_6,  2)$    \\ [0.8ex] 
 \hline
$ \fg^\sigma $  &  $ C_n $  &   $ C_n $   &    $B_n$ &  $B_n $   &    $G_2 $ &  $F_4 $  \\ [0.8ex] 
 \hline
\end{tabular},
\end{equation}
where by convention $C_1$ and $B_1$ are $A_1$ and $n\geq 3$ for $D_{n+1}$.

\begin{lemma}
\label{A_2n_lem}
Let $\fg$ be of type $A_{2n}$ and let $\sigma$ be the standard nontrivial automorphism of $\fg$. Then, the bijection $ D_{c, \sigma}\simeq  D_{c,\tau}$ in (\ref{weight_formula})
is given by $\sum_{i=1}^n  a_i\lambda_i^C\mapsto \sum_{i=1}^{n-1} a_i \lambda_i^{B}  + (2a_n+c )\lambda_{n}^B$, where $\{ \lambda_i^C\,|\, i=1,\cdots, n  \}$ is the set of fundamental weights of $\fg^\sigma$, and $\{ \lambda_i^B\,|\, i=1,\cdots, n  \}$ is the set of fundamental weights of $\fg^\tau$. (We follow the labelings in \cite[Table Fin, p. 53]{Ka}.)
\end{lemma}
\begin{proof}
Let $\alpha_1,\cdots, \alpha_{2n}$ be the set of simple roots of $\fg$, and let $\check{\alpha}_1,\cdots, \check{\alpha}_{2n}$ be the set of simple coroots of $\fg$. Then, $\alpha_1|_{\fh^\sigma},\cdots, \alpha_{n-1}|_{\fh^\sigma}, 2\alpha_n|_{\fh^\sigma}$ form a set of simple roots for $\fg^\sigma$, cf.\,\cite[$\S$ 2.1]{BH}, and $\{ \check{\alpha}_i+\check{\alpha}_{2n+1-i} \,|\, i=1, \cdots, n   \}$ is the set of simple coroots of $\fg^\sigma$. On the other hand, $\{ \alpha_1|_{\fh^\sigma},\cdots, \alpha_{n-1}|_{\fh^\sigma}, \alpha_n|_{\fh^\sigma}\}$ form a set of simple roots for $\fg^\tau$, and $\{ \check{\alpha}_1+\check{\alpha}_{2n}, \cdots, \check{\alpha}_{n-1}+\check{\alpha}_{n+2}, 2(\check{\alpha}_n+\check{\alpha}_{n+1} )\}$ form the set of simple coroots of $\fg^\tau$. The lemma now easily follows from \cite[Formula (6)]{HK}.
\end{proof}

\subsection{Affine Weyl group of twisted affine Lie algebras}
\label{aff_weyl_section}
Let $\sigma$ be a standard nontrivial automorphism of $\fg$ and 
 let $\tilde{L}(\fg,\sigma)$ be the Lie algebra $\hat{L}(\fg,\sigma)\oplus \mathbb{C}d$, where 
\[ [d, x[t^k]]=k x[t^k],  \quad [d, C]=0 ,     \quad \text{ for any }    x[t^k]\in \hat{L}(\fg,\sigma).   \]
Then, $\tilde{L}(\fg,\sigma)$ is a  Kac-Moody Lie algebra of twisted type with canonical center $C$ and the scaling element $d$, and the fixed  subalgebra $\fg^\sigma$ is the ``standard" finite part of $\tilde{L}(\fg,\sigma)$  in the sense of \cite[\S 6.3]{Ka}. This is obvious when $\fg$ is not of type $A_{2n}$. When $\fg$ is of type $A_{2n}$, using the formula (\ref{as_auto})  this can be seen from \cite[Theorem 8.7]{Ka} or \cite[$\S$ 2]{HK}.  

Set $\tilde{\fh}:= \fh^\sigma \oplus \bbC C \oplus \bbC d$. Then, the dual $\tilde{\fh}^*= (\fh^\sigma)^*\oplus  \bbC \delta \oplus  \bbC \Lambda_0$, where $\delta $ and $\Lambda_0$ are defined as follows
\[ \delta|_{\fh^\sigma}=0,  \, ( \delta,C )=0,\   (  \delta, d  )=a_0,    \quad 
 \Lambda_0|_{\fh^\sigma}=0, \,   ( \Lambda_0, C)=1, \, ( \Lambda_0, d )=0,  
 \]
where
\[ a_0=\begin{cases}  
1 \quad  \text{ if } (\fg, m)\not=(A_{2n}, 4)\\
2  \quad   \text{ if } (\fg, m)=(A_{2n}, 4).
 \end{cases}\]  
 Note that $\Lambda_0$ is a fundamental dominant weight of $\tilde{L}(\fg,\sigma)$ of level one.   
 

Let $W_{\tilde{L}(\fg,\sigma)  }$ denote the Weyl group of $\tilde{L}(\fg,\sigma)$.  Let $Q_\sigma$ (resp. $P_\sigma$) be the root (resp. weight) lattice of $\fg^\sigma$.  Set 
\begin{equation}
\label{lattice_M}
M=\begin{cases}
Q_\sigma    \quad     \text{ if }  (\fg, m)\not=(A_{2n}, 4)\\
\frac{1}{2} Q_{\sigma, l}  \quad  \text{ if }   (\fg, m)=(A_{2n}, 4),
\end{cases}
\end{equation}
where $Q_{\sigma, l}$ is the lattice spanned by the long roots.  Let $W$ be the Weyl group of $\fg$. Then, the fixed  subgroup $W^\sigma$ can be identified with the Weyl group of $\fg^\sigma$.   
Let $W_{\sigma,c}$ denote the affine Weyl group $W^\sigma \ltimes c M$. 


 Set 
 \[\fh^*_{\sigma,\bbR} := P_\sigma\otimes_\bbZ \bbR , \quad   \tilde{\fh}^*_\bbR:=P_\sigma\otimes_\bbZ \bbR+\bbR \Lambda_0+\bbR \delta.\]
  Note that $W_{\tilde{L}(\fg,\sigma)  }$  keeps $\delta$ invariant (cf.\,\cite[$\S$6.5]{Ka}). Hence, $W_{\tilde{L}(\fg,\sigma)  }$  acts on $ \hat{\fh}^*_{\bbR, c}$  
  for any $c\in \bbR$, where 
 \[\hat{\fh}^*_{\bbR, c} :=  \{ x\in \tilde{\fh}^*_\bbR  \, |\,   ( x, C  )=c   \} / \bbR \delta. \]

With respect to the isomorphism  $\fh^*_{\sigma,\bbR}  \simeq \hat{\fh}^*_{\bbR, c}$
 given by $\lambda\mapsto  c \Lambda_0+ \lambda$, we have the following lemma (cf.\,\cite[Prop. 6.5, \S 6.6]{Ka} or \cite[Lemma 3.1]{Ho1}).


 \ble
 \label{dot_lem}
 There exists an isomorphism ${\rm af}: W_{\tilde{L}(\fg,\sigma)  } \simeq  W_{\sigma,c}$ of groups such that 
 for any  $\Lambda=c \Lambda_0 +\lambda\in \hat{\fh}^*_{\bbR, c}$ with $\lambda \in \fh^*_{\sigma,\bbR}$
and $w\in W_{\tilde{L}(\fg,\sigma) } $, 
 the following formula holds,
 \[  w\cdot \Lambda=c \Lambda_0+  {\rm af}(w)\cdot  \lambda \text{ in } \hat{\fh}^*_{\bbR, c}. \]
 \ele
%
 
 Let $\hat{\rho}$ be the sum of all the fundamental weights of $\tilde{L}(\fg,\sigma)$. By \cite[\S 6.2.8]{Ka}, 
 \beq 
 \label{hatrho} \hat{\rho}=\rho_\sigma+\check{h}\Lambda_0,
 \eeq  where $\rho_\sigma$ is the sum of all the fundamental weights of $\fg^\sigma$, and $\check{h}$ is the dual Coxeter number of $\tilde{L}(\fg,\sigma)$, cf.  \,\cite[\S 6.1]{Ka}. Observe that $\check{h}$ is the same as  the dual Coxeter number of $\fg$ (cf. \cite[Remark 6.1]{Ka}).   
 
 We define $\star$ action of $W_{\tilde{L}(\fg,\sigma)  }$  on $ \hat{\fh}^*_{\bbR, c}$ as follows:
 \[ w\star \Lambda=w\cdot (\Lambda+\hat{\rho})-\hat{\rho}, \quad  w\in W_{\tilde{L}(\fg,\sigma)  } , \Lambda\in \hat{\fh}^*_{\bbR, c}. \]
Similarly, we still denote by  $\star$ the following action of $W_{\sigma, c} $ on $\fh^*_{\sigma,\bbR} $:
 \beq
 \label{star_act}
  w\star \lambda=w\cdot(\lambda+\rho_\sigma)-\rho_\sigma, \quad  w\in W_{\sigma,c}, \lambda\in \fh^*_{\sigma,\bbR}. \eeq
\ble
\label{star_leml}
Given $\Lambda=c \Lambda_0+\lambda\in \hat{\fh}^*_{\bbR, c} $ and $w\in W_{\tilde{L}(\fg,\sigma)  } $, we have 
\[ w\star \Lambda=\lambda+ {\rm af}_{\check{h}}(w)\star \lambda, \quad \text{ where }  {\rm af}_{\check{h}}(w)\, \,\text{is taken in}\,  W_{\sigma, c+\check{h}}   . \]
\ele
\bpf
It follows from Lemma \ref{dot_lem} and the formula  $\hat{\rho}=\rho_\sigma+\check{h}\Lambda_0$ as in \eqref{hatrho}.
\epf

Set
\beq
\label{theta_1}
 \theta_\sigma=\begin{cases}
\text{ highest short root of   } \fg^\sigma,    \quad  (\fg, m)\not= (A_{2n}, 4)\\
\frac{1}{2}\text{ highest  root of $\fg^\sigma$},    \quad  (\fg, m)= (A_{2n}, 4)

 \end{cases}   \eeq
and
\beq
\label{theta_2}
\check{ \theta}_\sigma=\begin{cases}
\text{ highest coroot of   } \fg^\sigma,    \quad  (\fg, m)\not= (A_{2n}, 4)\\
2\cdot \text{ highest short coroot of  } \fg^\sigma ,   \quad  \quad  (\fg, m)= (A_{2n}, 4).
 \end{cases}   \eeq
 Let $\langle\cdot |\cdot \rangle$ denote the normalized bilinear form on ${\fh}^\sigma$ (which is the restriction of the normalized invariant form on $\fh$). Let $\nu\colon \fh^\sigma\simeq (\fh^\sigma)^*$ be the induced isomorphism.  Then, $\nu(\check{\theta}_\sigma )= \frac{1}{a_0}\theta_\sigma $, cf. \,\cite[\S 6.4]{Ka}.  

When $\sigma$ is standard and nontrivial, using  formula (\ref{as_auto})  (for the case $ (\fg, m)= (A_{2n}, 4)$) combined with 
\cite[$\S$ 2]{HK} for the diagram automorphisms, we get:
\beq 
\label{ab_weight}
 D_{c,\sigma} =\{  \lambda\in P^+_\sigma  \,|\,    (  \lambda,   \check{\theta}_\sigma )  \leq  c \ra \}   ,\eeq
where $P_\sigma^+$ is the set of dominant integral weights of $\fg^\sigma$.

Let $ W_{\sigma, c+\check{h} }^\dagger$ denote the set of minimal representatives of the left cosets of $W^\sigma$ in $W_{\sigma, c+\check{h} }$. Then, for any $w_1\in W^\sigma$ and $w_2\in W_{\sigma, c+\check{h} }^\dagger $, we have $\ell(w_1w_2)=\ell(w_1)+\ell(w_2)$. For any $w\in W^\dagger_{\sigma, c+ \check{h}}$ and $\lambda\in D_{c,\sigma}$, one has $w\star \lambda\in P_\sigma^+$, and $w\star \lambda\not= w'\star \lambda$ for any two distinct $w,w'\in W^\dagger_{\sigma, c+ \check{h}}$ (cf.\,\cite[Remark 1.3]{Ko}).  Since  $( \rho_\sigma,  \check{\theta}_\sigma )= \check{h}-1$, $D_{c,\sigma}$ can be identified with the interior integral points in the fundamental alcove of $W_{\sigma, c+\check{h}}$ with respect to the $\star$ action.



\subsection{Analogue of Kac-Walton formula}
\label{subsect_KW}

Let $\Gamma=\la \sigma \ra$ of order $m$ act on $\mathbb{P}^1$  by $\sigma(z)=e^{\frac{2\pi \mathrm{i}}{m}} z$ (for $z\in \mathbb{P}^1$), and $\sigma$ acts on $\fg$ via a standard automorphism of order $m$.  For any $z\in \bbP^1\backslash \{0\}$ and any finite dimensional representation $V$ of $\fg^{\Gamma_z}$, where $\Gamma_z$ is the stabilizer subgroup of $\Gamma$ at $z$, we denote by $V_z$ the representation of $\fg[t^{-1}]^\sigma$ via the evaluation map ${\rm ev}_z\colon \fg[t^{-1}]^\sigma\to \fg^{\Gamma_z}$ by letting $t=z$. 
 Recall that for any $\lambda\in D_{c,\sigma}$, we have an integrable highest weight representation $\mathscr{H}_c(\lambda)$ of $\hat{L}(\fg, \sigma)$ of level $c$ and highest weight $\lambda$. 
 Let $H_i( (t^{-1}\fg[t^{-1}])^\sigma,  \mathscr{H}_c(\lambda) \otimes V(\mu)_1)$ denote the $i$-th homology of  $(t^{-1}\fg[t^{-1}])^\sigma$
with coefficients in  $\mathscr{H}_c(\lambda) \otimes V(\mu) _1$. 
\bpr
\label{BGG_comples}
For any $\lambda\in D_{c,\sigma}$ and $\mu\in P^+$ (where $P^+$ is the set of dominant integral weights of $\fg$),  the homology groups  $H_*( (t^{-1}\fg[t^{-1}])^\sigma,  \mathscr{H}_c(\lambda) \otimes V(\mu) _1)$  can be computed as the homology groups of a complex of $\fg^\sigma$-representations,
   \[ \cdots   \to   F_p  \xrightarrow{\bar\delta_p}  \cdots   F_1   \xrightarrow{\bar\delta_1} F_0  \xrightarrow{\bar\delta_0}   0, \]
where as representations of $\fg^\sigma$,  
\beq 
\label{D_complex}
F_p\simeq  \bigoplus_{ w\in\,  W^\dagger_{\sigma, c+\check{h}} ,  \ell(w)=p   }    V(w\star \lambda) \otimes  (V(\mu) |_{\fg^\sigma}) ,  \eeq
 $V(w\star \lambda)$ is the irreducible representation of $\fg^\sigma$ with highest weight $w\star \lambda$ and $V(\mu)|_{\fg^\sigma}$ denotes   the irreducible representation $V(\mu)$ of $\fg$ considered as a representation of  $\fg^\sigma$ via restriction. 
 
\epr
\bpf

Recall the generalized BGG resolution for  Kac-Moody algebras from \cite[Definition 9.2.17]{Kbook}. By Lemmas \ref{dot_lem} and \ref{star_leml}, we can express the  generalized BGG resolution of $\mathscr{H}_c(\lambda)$ as follows: 
 \[ \cdots   \to   M_p   \xrightarrow{\delta_{p}}  \cdots   \xrightarrow{\delta_1}  M_0  \xrightarrow{\delta_0} \mathscr{H}_c(\lambda), \]
 where 
  \beq
\label{F_lambda}
 M_{p}:=   \bigoplus_{w\in\, W^\dagger_{\sigma, c+\check{h}}, \ell(w)=p} \hat{M}({w\star \lambda}),
 \eeq
 and $ \hat{M}({w\star \lambda})$ is the generalized Verma module $U( \hat{L}(\fg,\sigma)  )\otimes_{U ( \fg[t]^\sigma \oplus \bbC C  )  }  V(w\star \lambda)  $, with $C$ acting on $V(w\star \lambda) $  by the scalar $c$ and $\fg[t]^\sigma$ acting via the evaluation at $t=0$.

By tensoring with $V(\mu) _1$, we get a resolution of $\mathscr{H}_c(\lambda) \otimes V(\mu) _1$.  As $\fg^\sigma$-modules, we have the coinvariant
 \begin{align*}
  (\hat{M}({w\star \lambda }  ) \otimes V(\mu) _1)_{    (t^{-1}\fg[t^{-1}])^\sigma     }
  &\simeq \left( (U(  (t^{-1}\fg[t^{-1}])^\sigma)  \otimes_\bbC  V(w\star \lambda )
    )\otimes V(\mu) _1\right)_{ (t^{-1}\fg[t^{-1}])^\sigma   }
       \\
  &\simeq \left( ( V(w\star \lambda )\otimes V(\mu)_1)\otimes_\bbC U(  (t^{-1}\fg[t^{-1}])^\sigma    )  \right)_{ (t^{-1}\fg[t^{-1}])^\sigma   },\\
  &\medskip\medskip
  \,\,\,\text{by the Hopf Principle \cite[Proposition 3.1.10]{Kbook}}\\
  &\simeq  V(w\star \lambda )\otimes V(\mu).  
  \end{align*}
Hence, the complex $(M_\bullet \otimes V(\mu) _1)_{  (t^{-1}\fg[t^{-1}])^\sigma }$ is isomorphic to 
\[ \cdots   \to   F_p  \xrightarrow{\bar\delta_p}  \cdots   F_1   \xrightarrow{\bar\delta_1} F_0  \xrightarrow{\bar\delta_0}   0,\]
where $F_p$ is given in (\ref{D_complex}). Thus, the proposition follows. 
\epf

Consider the automorphism $\sigma$ of $\bbP^1$ given by $\sigma(z)=\xi z$ for $z\in \mathbb{P}^1$, where $\xi=e^{ \frac{2\pi \mathrm{i}}{m}}$. With respect to the Galois cover $\pi: \bbP^1\to \bbP^1$ given by $z\mapsto z^m$, the ramification type at $0\in \bbP^1$ is $\sigma$ and the ramification type at $\infty \in \bbP^1$ is $\sigma^{-1}$. For twisted conformal blocks associated to the Galois cover $\pi$ and the automorphism $\sigma$ of $\fg$, 
 we attach a dominant weight in $D_{c,\sigma}$  at $0\in \bbP^1$ , and  a dominant weight in $D_{c,\sigma^{-1}}$ at $\infty\in \bbP^1$.  We have the following lemma.
\ble
\label{duality_lem}
For any standard nontrivial automorphism $\sigma$ of $\fg$,  we have $D_{c,\sigma}=D_{c, \sigma^{-1}}$. Moreover, for any $\mu\in D_{c,\sigma}$, $\mu=\mu^*$ where $\mu^*$ is the dominant weight corresponding to the dual representation $V(\mu)^*$ of $\fg^\sigma$.
\ele
\bpf
By  \cite[Lemma 5.3 (2)]{HK}, $\mu\in D_{c,\sigma}$ if and only if $\mu^*\in D_{c,\sigma^{-1}}$.  Since $\fg^\sigma$ is non simply-laced or $A_1$,
$\lambda^*=-w_0^\sigma(\lambda)=\lambda$ for any weight $\lambda$ of $\fg^\sigma$, where $w_0^\sigma$ is the longest element of the Weyl group of $\fg^\sigma$. 
\epf
Thus, we can use $D_{c,\sigma}$ for the common set of dominant weights of $\fg^\sigma$ to attach to $0$ as well as  $\infty\in \bbP^1$.

Similar to Teleman's vanishing theorem  \cite[Theorem 0]{T1}, we make the following conjecture. 
\begin{conjecture}
\label{Teleman} Let $\sigma$ be a standard automorphism of $\fg$. 
Then, for any $\lambda, \mu \in    D_{c,\sigma}$, $ \nu\in D_c$,  and for any $i\geq 1$,  the representation $V(\mu)^*$ does not occur in $H_i((t^{-1}\fg[t^{-1}])^\sigma, \mathscr{H}_c(\lambda)\otimes  V(\nu) _1  )$ as a $\fg^\sigma$-representation.
\end{conjecture}

We are now ready to deduce the following analogue of Kac-Walton formula for twisted conformal blocks.
\bt
\label{KW_thm} Take any standard representation $\sigma$ of $\fg$ of order $m$ and the Galois cover $\pi: \bbP^1\to 
\bbP^1, z\mapsto z^m$. 
Let $\vec{p}=(0,\infty, 1) $ in $\bbP^1$, and $\vec{\lambda}=(\lambda, \mu, \nu)$ with $\lambda, \mu\in D_{c,\sigma}$ and $\nu\in D_c$.  Suppose that Conjecture \ref{Teleman} holds, 
then
\beq  \dim  \mathscr{V}_{\bbP^1, \Gamma, \phi}(\vec{p},  \vec{\lambda})=  \sum_{ w\in  W^\dagger_{\sigma, c+\check{h}}  }   (-1)^{\ell(w)  }   \dim \left((V(w\star \lambda)\otimes V(\mu)\otimes V(\nu) )^{\fg^\sigma} \right),        \eeq
where $\Gamma$ is the cyclic group of order $m$ generated by $\sigma$ and $\phi: \Gamma \to \Aut \fg$ is the one generated by $\sigma$.  
\et

\bpf
By the Propagation Theorem for twisted conformal blocks (cf.\,\cite[Theorem 4.3]{HK}), we have the following isomorphism:
\begin{align}\label{eqn4.7.1}
 \mathscr{V}_{\bbP^1, \Gamma, \phi}(\vec{p},  \vec{\lambda})
&\simeq   (\mathscr{H}_c(\lambda) \otimes V(\mu)_\infty\otimes  V(\nu) _1 )_{\fg[t^{-1}]^\sigma}   \notag  \\  
&\simeq  \left( (\mathscr{H}_c(\lambda) \otimes  V(\nu) _1 )_{(t^{-1}\fg[t^{-1}])^\sigma} \otimes V(\mu)\right)_{\fg^\sigma},\,\,\,\text{since   $(t^{-1}\fg[t^{-1}])^\sigma$ acts trivially on $V(\mu)_\infty$}\notag\\
&\simeq \Hom_{\fg^\sigma}\left(V (\mu^*), H_0((t^{-1}\fg[t^{-1}])^\sigma, \mathscr{H}_c(\lambda)\otimes V(\nu) _1 )\right).
\end{align}
We have the following equalities:
\begin{align}\label{eqn4.7.2}
  & \dim   \left(\Hom_{\fg^\sigma}(V (\mu^*), H_0((t^{-1}\fg[t^{-1}])^\sigma, \mathscr{H}_c(\lambda)\otimes V(\nu) _1 ))\right)\notag\\
   &=\sum_{i\geq 0} (-1)^i \dim \left(\Hom_{\fg^\sigma}(V(\mu^*),H_i((t^{-1}\fg[t^{-1}])^\sigma, \mathscr{H}_c(\lambda)\otimes V(\nu) _1 ) ) \right)\notag\\
   &=\sum_{w\in W^\dagger_{\sigma, c+ \check{h}} }   (-1)^{\ell(w)  }  \dim\left(  \Hom_{\fg^\sigma}(V(\mu^*), V(w\star \lambda) \otimes  V(\nu)) \right)\notag \\
   &=\sum_{w\in W^\dagger_{\sigma, c+ \check{h}}  }  (-1)^{\ell(w)  } \dim \left( (V(\mu)\otimes V(w\star \lambda) \otimes  V(\nu))^{\fg^\sigma} \right),
\end{align}
where the first equality follows from Conjectutre \ref{Teleman}, and the second equality follows from Proposition \ref{BGG_comples}. Combining the isomorphism \eqref{eqn4.7.1} and the identity \eqref{eqn4.7.2}, we get the theorem. 
\epf

\section{Verlinde formula for twisted conformal blocks}
\label{twisted_Verlinde_sect}

\subsection{Verlinde formula for basic cases}
In this section, we assume that $\sigma$ is a standard nontrivial automorphism of $\fg$.  Recall the lattice $M$ introduced in (\ref{lattice_M}).  Then, $M$ is the root lattice of $\fg^\sigma$, if $\fg$ is not $A_{2n}$, and by \cite[Lemma 2.3]{Ho1}, $M$ is the weight lattice of $\fg^\sigma$ when $\fg$ is of type $A_{2n}$. Let $G$ be the connected and  simply-connected (simple) algebraic group with Lie algebra $\fg$. Let $T$ be the maximal torus with Cartan subalgebra $\fh$ as its Lie algebra. 
\begin{lemma}\label{lemma5.1}
The fixed group $G^\sigma$ is connected and simply-connected. 
\end{lemma}
\begin{proof}
When $G$ is not of type $A_{2n}$, $\sigma$ is a diagram automorphism. In this case, the lemma is well-known.  We now assume that $G$ is of type $A_{2n}$.  Let $\tau$ be the diagram automorphism part of $\sigma$. Then, $T^\sigma=T^\tau$. It is known that $T^\tau$ is connected. Thus, $G^\sigma$ is connected of type $C_n$, see Table (\ref{Fix_table}). Let $\{\alpha_1,\cdots, \alpha_{2n}\}$ be the set of simple roots of $A_{2n}$ with the standard labelling and let  $\{\check{\alpha}_1,\cdots, \check{\alpha}_{2n}\}$ be the set of corresponding simple coroots of $G$. 
Then, $\alpha_1|_{\fh^\sigma},\cdots, \alpha_{n-1}|_{\fh^\sigma}, {2\alpha_n}_{|\fh^\sigma}$ form a set of simple roots for $G^\sigma$  (cf.\,\cite[$\S$ 2.1]{BH}) and $\{\check\alpha_i +\check\alpha_{2n+1-i}: 1\leq i\leq n\}$ form the set of simple coroots. Using simple coroots, we can introduce a coordinate system of $T$, $(\mathbb{G}_m)^{2n} \simeq T$. Then, it is easy to verify that $\mathbb{G}_m\to T^\sigma$ given by $a\mapsto \check{\alpha}_{i}(a)\check{\alpha}_{2n+1-i}(a)$ is a simple coroot of $G^\sigma$, for every $i=1,\cdots, n$. Thus, ${\rm Hom}(\mathbb{G}_m, T^\sigma )$ is the lattice of coroots.  It follows that $G^\sigma$ is simply-connected. 
\end{proof}

With this lemma, we may regard the lattice $M$ as a sub-lattice in the weight lattice $X^*(T^\sigma)$. We now define 
\begin{equation}
\label{sigma_torus_set}
  T^\sigma_c: =\{ t\in T^\sigma \,|\,   e^{\lambda} (t)=1, \, \forall  \lambda \in     (c+ \check{h}) M    \}   .  \end{equation}

Let $T^{\sigma, reg}_c$ denote the set of regular elements in $T^\sigma_c$, i.e.,  those elements with trivial $W^\sigma$-stabilizer.  Let $R^\sigma_c(\fg )$ denote the fusion ring associated to the twisted affine Lie algebra $\hat{L}(\fg,\sigma)$, which is defined in \cite{Ho2}. For any regular function $f$ on $T^\sigma$ or $T$, we will denote by $\bar{f}$  the restriction of $f$ to $T^{\sigma, reg}_c$. 
 In \cite{Ho2}, the ring $R^\sigma_c(\fg)$ is realized as the function space $\bbC[T_c^{\sigma, reg}/W^\sigma]=\mathbb{C}[T^{\sigma, reg}_c]^{W^\sigma}$ (with the ring structure coming from the product of functions), with a basis $\{ \bar{\chi}_\lambda \,|\, \lambda\in D_{c,\sigma} \}$. (We describe  $\bar{\chi}_\lambda$ explicitly after Remark  \ref{remark6.3}.)
 The following theorem is proved in \cite{Ho1,Ho2}. 
 \begin{theorem}
  \label{verlinde_Hong}
For any $\lambda, \mu\in D_{c,\sigma}$, we have 
 \begin{equation*} 
  \bar{\chi}_\lambda \cdot \bar{\chi}_\mu= \sum_{\eta \in D_{c,\sigma} } c_{\lambda,\mu}^\eta \bar{\chi}_\eta ,   \end{equation*}
where $c_{\lambda, \mu}^\eta$ is given by 
\begin{equation}\label{verlinde_eq_h}
c_{\lambda, \mu}^\eta=  \frac{1}{| T^\sigma_c |} \sum_{t\in T^{\sigma, reg}_c/W^\sigma  }   \bar\chi_\lambda(t) \bar\chi_\mu(t) \bar\chi_{\eta^* }(t)  \Delta_\sigma(t) . \end{equation}
Here $\Delta_\sigma$ is given by 
\[\Delta_\sigma := \prod_{\alpha\in \Phi_\sigma}  (e^\alpha-1),  \]
where $\Phi_\sigma$ is the set of all the roots of $\fg^\sigma$. 

 \end{theorem}

\begin{remark}\label{remark6.3}
{\rm Given a simply-laced simple Lie algebra $\dot{\fg}$ with a diagram automorphism $\dot{\tau}$ of order $\dot r >1 $, a fusion ring $R_c(\dot{\fg}, \dot{\tau} )$ is defined in \cite{Ho1} for the purpose of deducing a  formula for the trace of $\dot{\tau}$ on the space of untwisted conformal blocks associated to $\dot{\fg}$. In fact, there is an isomorphism of rings $R_c^\sigma(\fg)\simeq R_c(\dot{\fg}, \dot{\tau} ) $, with a correspondence between $(\fg,m)$ and $(\dot{\fg},\dot{r})$ as follows (cf.\,\cite[$\S$ 3.1]{Ho2}):
 \begin{equation}
\label{corresp_table}
 \begin{tabular}{|c  | c | c |c |c |c | c| c| c|c |c|c|c|c|c |c ||} 
 \hline
$(\fg, m)$   & $(A_{2n-1}, 2 ) $  &  $(A_{2n}, 4)$  &  $(D_{n+1} , 2  ) $  &           $ (D_4,   3)$  &  $ (E_6,  2)$    \\ [0.8ex] 
 \hline
$ (\dot{\fg}, \dot{r}) $  &  $ (D_{n+1},2) $  &   $ (A_{2n}, 2) $  &  $ (A_{2n-1}, 2 )$   &    $(D_4, 3) $ &  $(E_6, 2) $  \\ [0.8ex] 
 \hline
\end{tabular}.
\end{equation}
Moreover, Theorem \ref{verlinde_Hong} is equivalent to the formula for the trace of $\dot{\tau}$ on the space of conformal blocks associated to $\dot{\fg}$. }
\end{remark}

There is a natural ring homomorphism $\pi\colon R(\fg^\sigma)\to R^\sigma_c(\fg)$ given by $[V_\lambda]\mapsto \bar{\chi}_\lambda$, where $R(\fg^\sigma)$ is the representation ring of $\fg^\sigma$,  $\chi_\lambda$ is the character of $V(\lambda)$ as a function on $T^\sigma$ and $\bar\chi_\lambda$ is the restriction of $\chi_\lambda$ to $T_c^{\sigma, reg}$ (which descends to a function on $T_c^{\sigma, reg}/W^\sigma$). 
 
\begin{lemma}
\label{sign_lem}
For any $\lambda\in D_{c,\sigma}$ and $w\in\, W^\dagger_{\sigma, c+\check{h}}$ ,  we have
\[  \bar{\chi}_{w\star \lambda}  = (-1)^{\ell(w)} \bar{\chi}_{\lambda } .    \]
Moreover, by \cite[Corollary 5.17]{Ho1}, for any $\eta \in P_\sigma^+\setminus (W^\dagger_{\sigma, c+\check{h}}\star  D_{c,\sigma})$, $\pi(V_\eta)=0.$
\end{lemma}
\begin{proof}
Write $w=z\tau_{\eta }$, where $z\in W^\sigma$ and $\tau_\eta$ is the translation by $\eta\in (c+\check{h})M $. 
By Weyl character formula, 
\begin{equation}
\delta_\sigma \cdot  \chi_{w\star \lambda} = \sum_{y\in W^\sigma}  (-1)^{\ell(y)} e^{ y(w\star \lambda+\rho_\sigma) } =\sum_{y\in W^\sigma}  (-1)^{\ell(y)} e^{ y(z( \lambda+\rho_\sigma+\eta ) ) },
\end{equation}
where $\delta_\sigma$ is the Weyl denominator of $\fg^\sigma$ given by:
\[\delta_\sigma := e^{\rho_\sigma} \prod_{\alpha\in \Phi_\sigma^+}  (1-e^{-\alpha})   \]
 ($\Phi_\sigma^+$ being the set of the positive  roots of $\fg^\sigma$).

For any $t\in T^{\sigma, reg}_c$, we have 
\begin{align*}   \delta_\sigma(t) \cdot  \chi_{w\star \lambda}(t) &=  \sum_{y\in W^\sigma}  (-1)^{\ell(y)} e^{ y(z( \lambda+\rho_\sigma ) ) }(t)=(-1)^{\ell(z)} \sum_{y\in W^\sigma}  (-1)^{\ell(y)} e^{ y(\lambda+\rho_\sigma) } (t) \\
& =(-1)^{\ell(w)} \sum_{y\in W^\sigma}  (-1)^{\ell(y)} e^{ y(\lambda+\rho_\sigma) }(t)= (-1)^{\ell(w)}\delta_\sigma(t) \cdot \chi_{ \lambda}(t),
\end{align*}
where the first equality holds since $\eta\in (c+\check{h})M$, and 
second to the last  equality holds since $\ell(\tau_\eta)$ is even, cf.\,\cite[Lemma 2.8]{Ho1}. Thus, the lemma follows. 
  
\end{proof}

Let $\vec{p}=(0,\infty, 1) $ in $\bbP^1$, and $\vec{\lambda}=(\lambda, \mu, \nu)$ with $\lambda, \mu\in D_{c,\sigma}$ and $\nu\in D_c$.  Recall the following notation from Notation \ref{verlinde_notation} (we have dropped $\phi$ from 
$N_\phi(\sigma; \lambda,\mu,\nu)$ since in this section we are only dealing with $\phi$ generated by the nontrivial standard automorphisms of $\fg$): 
\begin{equation}
\label{dim_verlinde}
N(\sigma; \lambda,\mu,\nu)= \dim   \mathscr{V}_{\bbP^1, \Gamma, \phi}(\vec{p},  \vec{\lambda}), 
\end{equation}
where $\Gamma$ and $\phi$ are the same as in Theorem \ref{KW_thm}.

We now prove the following Verlinde formula for $N(\sigma; \lambda,\mu,\nu)$, which uses Theorem \ref{KW_thm}. 
\bt
\label{verlinde_genus_0}
With the notation as above, suppose that Conjecture \ref{Teleman} holds. Then, we have
\beq   N(\sigma;\lambda,\mu,\nu)= \frac{1}{| T^\sigma_c |} \sum_{t\in T^{\sigma, reg}_c/W^\sigma  }    \chi_\lambda(t)\chi_\mu(t) \chi_\nu(t) \Delta_\sigma(t),      \eeq
where $\chi_\lambda, \chi_\mu, \chi_\nu$ represent the characters of $V(\lambda), V(\mu), V(\nu)$ as representations of $\fg^\sigma$, $\fg^\sigma$
 and $\fg$ respectively. 
\et

\bpf

For any $\mu\in D_{c,\sigma}$ and $\nu\in D_c$, 
consider the following decomposition as representations of $\fg^\sigma$,
\[  V(\mu)\otimes (V(\nu)|_{\fg^\sigma}) = \bigoplus_{ \eta \in P^+_\sigma}  V(\eta)^{\oplus  m_{ \mu,\nu} ^\eta } ,\]
where $m_{ \mu,\nu} ^\eta := \dim \left({\rm Hom}_{\fg^\sigma}(V(\eta), V(\mu)\otimes V(\nu) )\right)$.  Then, 
\[ V(\mu)\otimes (V(\nu)|_{\fg^\sigma}) = \left(\bigoplus_{\lambda \in D_{c,\sigma}}  \bigoplus_{w\in\, W^\dagger_{\sigma, c+\check{h}}    }   V(w\star \lambda) ^{\oplus  m_{\mu,\nu}^{w\star \lambda }  } \right)\bigoplus \left(\bigoplus_{\eta \in P_\sigma^+\setminus (W^\dagger_{\sigma, c+\check{h}}\star  D_{c,\sigma})}\,
V(\eta)^{\oplus  m_{\mu,\nu}^{\eta} }\right).
  \]
Thus, by Lemma \ref{sign_lem}, we have 
\beq\label{eqn5.5.1}
 \pi( V(\mu)\otimes (V(\nu)|_{\fg^\sigma}) )= \sum_{\lambda \in D_{c,\sigma}}  \sum_{w\in\, W^\dagger_{\sigma, c+\check{h}}    }   m_{\mu, \nu}^{w\star \lambda }  \bar{\chi}_{ w\star \lambda }   = \sum_{\lambda\in D_{c,\sigma}}  \sum_{w\in\, W^\dagger_{\sigma, c+\check{h}}    }    (-1)^{\ell(w)}   m_{\mu, \nu}^{w\star \lambda }  \bar{\chi}_{\lambda }   . 
\eeq

By the analogue of the Kac-Walton formula as  in Theorem \ref{KW_thm},  we get using the equation \eqref{eqn5.5.1}:
\begin{equation}
\label{com_formula1}
 \pi( V(\mu)\otimes (V(\nu)|_{\fg^\sigma}) )=  \sum_{ \lambda\in D_{c,\sigma} }  N(\sigma;\lambda^*, \mu, \nu)  \bar{\chi}_\lambda  .
 \end{equation}
(Observe that, by Lemma \ref{duality_lem} $\lambda^*=\lambda$.)

On the other hand, consider the following decomposition as $\fg^\sigma$-representations: 
\begin{equation}
\label{branch_dec}
  V_\nu|_{\fg^\sigma}=  \bigoplus_{ \eta\in P_\sigma^+  }   V_{\eta}^{\oplus  b_{\nu }^\eta  } =\left( \bigoplus_{\lambda' \in D_{c,\sigma} }  \bigoplus_{w\in\, W^\dagger_{\sigma, c+\check{h}}    }  V(w\star \lambda') ^{\oplus  b_{\nu}^{w\star  \lambda' } }\right)  \bigoplus \left(\bigoplus_{\eta \in P_\sigma^+\setminus (W^\dagger_{\sigma, c+\check{h}}\star  D_{c,\sigma})}\,
V(\eta)^{\oplus b_\nu^{\eta} }\right).  
 \end{equation}

Then, $\pi$ being a ring homomorphism,
\beq
\label{eqn5.5.2}  \pi( V(\mu)\otimes (V(\nu)|_{\fg^\sigma}) )=  \sum_{\lambda' \in D_{c,\sigma} }  \sum_{w\in\, W^\dagger_{\sigma, c+\check{h}}    }  b_{\nu}^{w\star  \lambda' }     (-1)^{\ell(w)} \bar{\chi}_\mu \cdot \bar{\chi}_{\lambda'} .   \eeq
By  Theorem \ref{verlinde_Hong} and equation \eqref{eqn5.5.2}, $ \pi( V(\mu)\otimes (V(\nu)|_{\fg^\sigma}) )$ is equal to 

\begin{align}
&  \sum_{\lambda' \in D_{c,\sigma} }  \sum_{w\in\, W^\dagger_{\sigma, c+\check{h}}    }  b_{\nu}^{w\star  \lambda' }     (-1)^{\ell(w)}   \sum_{\lambda \in D_{c,\sigma}} \big ( \frac{1}{| T^\sigma_c |} \sum_{t\in T^{\sigma, reg}_c/W^\sigma  } \chi_{\lambda^*}(t)  \chi_\mu(t) \chi_{\lambda'}(t)  \Delta_\sigma(t) \big ) \bar{\chi}_\lambda \nonumber \\
 &=   \sum_{\lambda\in D_{c,\sigma} }  \frac{1}{| T^\sigma_c |} \sum_{t\in T^{\sigma, reg}_c/W^\sigma  }   \chi_{\lambda^*}(t) \chi_{\mu}(t) \chi_{\nu }(t)  \Delta_\sigma(t)  \bar{\chi}_\lambda  ,\label{com_formula2}
  \end{align}
  where the above equality follows from (\ref{branch_dec})  
   and Lemma \ref{sign_lem}.
Comparing formulae (\ref{com_formula1}) and (\ref{com_formula2}), we conclude that 
\[  N(\sigma; \lambda^*,\mu,\nu)= \frac{1}{| T^\sigma_c |} \sum_{t\in T^{\sigma, reg}_c/W^\sigma  }    \chi_{\lambda^*}(t)\chi_\mu(t) \chi_\nu(t) \Delta_\sigma(t) .\]

Thus, the theorem follows. 
\epf

Following \cite{Ho2}, we now describe the set $T^{\sigma, reg}_c/W^\sigma $ explicitly.  Let $\theta_l$ be the highest root of $\fg^\sigma$.  Let $ \check{P}^+_\sigma$ denote the set of dominant coweights of $\fg^\sigma$, where the fundamental coweights in $\fh^\sigma$ are defined as the dual of simple roots.  
When $(\fg, m)\not=(A_{2n}, 4)$, set 
\[  \check{D}_{c,\sigma}=  \{  \check{\lambda}\in  \check{P}^+_\sigma \, | \,  ( \check{\lambda}, \theta_l )  \leq c     \} ,  \]
and 
\[  \Sigma_c:= \{  e^{ \frac{2\pi i}{ c+\check{h}} ( \check{\rho}_\sigma+\check{\lambda}, \,\cdot \, ) } \in  T^\sigma \,|\,  \check{\lambda} \in \check{D}_{c,\sigma}      \} \subset  T^{\sigma, reg}_c .\]
 Here we identify $T^\sigma =\Hom(P_\sigma, \mathbb{C}^*)$,  where $P_\sigma$ is the weight lattice of $\fg^\sigma$
 and $ \check{\rho}_\sigma$ is the sum of the fundamental coweights.
 
When  $(\fg, m)=(A_{2n}, 4)$, set 
\[\Sigma_c: = \{  e^{ \frac{2\pi i}{ c+\check{h}} \langle \rho_\sigma+\lambda | \,\cdot \,  \rangle } \in  T^\sigma \,|\,  \lambda \in D_{c,\sigma}      \}  \subset T^{\sigma, reg}_c ,\]
where $\langle\cdot | \cdot \rangle$ is the invariant  form on $(\fh^\sigma)^*$ such that $\langle\theta_l | \theta_l\rangle=4$, equivalently $\langle\cdot | \cdot \rangle$ is induced from the normalized invariant form on the twisted affine algebra $\hat{L}(\fg,\sigma)$ (cf. \cite[Identity8.3.8]{Ka}). 
The following lemma follows from \cite[$\S$ 5.4]{Ho1} and \cite[$\S$ 2]{Ho2}. 
\begin{lemma}
Any element $t\in T^{\sigma, reg}_c$ can be translated to a unique element in $\Sigma_c$ by a unique element of $W^\sigma$. 
\end{lemma}
\subsection{Verlinde formula for general $\Gamma$-curves}

{\it Let $\sigma$ be a standard nontrivial automorphism of a simple and simply-connected algebraic group $G$ preserving a maximal torus $T$}.  Set 
\[ T_c=\{ t\in T   \,|\,  \lambda(t)=1,  \text{ for any } \lambda\in (c+\check{h})Q_{lg}  \, \},   \]
where $Q_{lg}$ is the sublattice of the root lattice of $G$ generated by the long roots (if all the root lengths are equal, we call them long roots).  Let $T_c^{reg}$ be the set of regular elements in $T_c$, i.e.,  those elements $t\in T_c$ whose stabilizers in the Weyl group $W$ is trivial. 
Recall that $T_{c}^\sigma$ is defined in (\ref{sigma_torus_set}). 
\begin{lemma}
\label{lem_5.2}
\begin{enumerate}
\item  $T_{c}^\sigma$ is the set of $\sigma$-invariants in $T_c$. \\
\item   The set $T_{c}^{\sigma, reg}/ W^\sigma$ can be identified with the set of $\sigma$-invariants in $T^{reg}_{c}/W$. 
\end{enumerate}
\end{lemma}
\begin{proof}
For part (1), it suffices to check that the lattice $M$ defined in (\ref{lattice_M}) is exactly the set of coinvariants of $\sigma$ in $Q$. When $\fg$ is not of type $A_{2n}$, this is obvious. When $(\fg, m)=(A_{2n}, 4)$, this follows from the description of the simple roots of
$\fg^\sigma$ in \cite[$\S$ 2.1]{BH} and \cite[Lemma 2.2]{Ho1}.  

We now prove part (2).  Since $G$ and $G^\sigma$ are simply-connected (cf. Lemma \ref{lemma5.1}), there exist $\sigma$-equivariant bijections 
\[    {\check{P}}/(c+\check{h})\check{Q} \simeq T_c, \quad    M^\vee/ {(c+\check{h})\check{Q}^\sigma  }\simeq T _c^\sigma,     \]
given by $\lambda\mapsto e^{\frac{2\pi i}{ c+\check{h} }  \lambda }$, where $\check{P}, \check{Q}$ are respectively the coweight and coroot lattices of $G$,  $\check{Q}^\sigma$ is the coroot lattice of $G^\sigma$ (it can also be identified with the set of $\sigma$-invariants in $\check{Q}$), and $M^\vee\subset \fh^\sigma$ is the dual lattice of $M \subset (\mathfrak{h}^\sigma)^*$.  In particular, $M^\vee$ is the coweight lattice of $G^\sigma$ when $G$ is not of type $A_{2n}$;   $M^\vee$ is the coroot lattice of $G^\vee$ when $G$ is of type $A_{2n}$, as in this case $M$ is the weight lattice of $G^\sigma$.  From the descriptions of coroots and coweights of $G^\sigma$ in \cite[$\S$ 2.1]{BH}, we observe that in any case $M^\vee= (\check{P})^\sigma$. 

Then, $T_c^{reg}/ W$ can be identified with the set of interior $\check{P}$-integral points
in the fundamental alcove of the affine Weyl group $W\ltimes (c+\check{h})\check{Q}$ (cf. \cite[Lemma 4.2.6 (b)]{Kbook2}). Similarly, $T_c^{\sigma, reg}/W^\sigma$ can be identified with the set of interior $M^\vee$-integral points in the fundamental alcove of the affine Weyl group $W^\sigma\ltimes  (c+\check{h})\check{Q}^\sigma$. By the same proof as in \cite[Prop.2.7]{Ho1},  the natural map 
$T_c^{\sigma, reg}/W^\sigma\to T_c^{reg}/W$ is a bijection. 
\end{proof}

Given any two finite order automorphisms $\gamma, \gamma'$ of $\fg$ such that they have the same images in ${\rm Out} (\fg)$, we can naturally identify $D_{c, \gamma}$ and $D_{c, \gamma'}$. More precisely, we first decompose $\gamma= \tau\epsilon ^{ad h} $ with respect to a $\gamma$-stable pair $(\mathfrak{b},\mathfrak{h})$,  and decompose $\gamma'=\tau' \epsilon^{ad h'}$ with respect to a $\gamma'$-stable pair $(\mathfrak{b}', \mathfrak{h}')$,  as in (\ref{eq1.1.1.0}).  By (\ref{weight_formula}), there exists identifications $D_{c, \gamma}\simeq D_{c, \tau}$ and $D_{c, \gamma'}\simeq D_{c, \tau'}$.  Furthermore, by Lemma \ref{lem_kappa_innner}, there exists a canonical identification $D_{c, \tau}\simeq D_{c, \tau'}$ (note that this identification does not depend on the choice of the inner automorphisms in Lemma \ref{lem_kappa_innner}).  Thus, we get an identification 
\begin{equation}
\label{wt_bij}
 D_{c, \gamma}\simeq  D_{c, \gamma'}.\end{equation}


Consider a group homomorphism $\phi: \Gamma\to  {\rm Aut}(\fg)$, a stable smooth (and hence irreducible) $s$-pointed $\Gamma$-curve $(\Sigma,\vec{p})$ with ramification type $\vec{\gamma}$ attached to $\vec{p}$, and a $s$-tuple of dominant weights
$\vec{\lambda}$ attached to $\vec{p}$. We assume that 
\begin{assumption}
\label{assump_of_thm}
\begin{enumerate}
\item  $\Gamma\cdot \vec{p}$ contains all the ramified points;
\item  $\tilde{\Gamma} :=\Gamma/\Gamma_0$ is cyclic of order $r$, where $\Gamma_0$ is the kernel of the map $P\circ \phi: \Gamma \to \Out (\fg); P$ being the projection ${\rm Aut}
 (\fg) \to \Out (\fg)$;  
\item    By reordering $\vec{p}$, we can write $\vec{p}= (p_1,\cdots, p_{2a},p_{2a+1}, \cdots, p_s )$ such that $\tilde{\gamma}_k\not= 1$, $\tilde{\gamma}_{2k-1}\tilde{\gamma}_{2k}=1$ for any $1\leq k\leq a$ ($a \geq 0$), and $\tilde{\gamma}_i=1$ for any $2a+1\leq i\leq s$, where $\tilde{\gamma}$ is the image of $\gamma\in \Gamma$ in $\tilde{\Gamma}$. 
\end{enumerate}
\end{assumption}
When $\tilde{\Gamma}= (1)$ or $\mathbb{Z}/2\mathbb{Z}$, the condition (3) in the  above assumption holds automatically.

Let $\tilde{\phi}_\iota: \tilde{\Gamma}\to {\rm Aut}(\fg)$ be the group homomorphism defined in Section \ref{reduction}, which preserves a fixed  pair $(\fb,\fh)$. For any $1\leq i\leq s$, $\tilde{\gamma}_i$ acts on $G$ by a diagram automorphism (possibly trivial) via $\tilde{\phi}_{\iota}$, denoted by $\tau_i$.  Let $\sigma_i$ be the standard automorphism associated to $\tau_i$ (note that $\sigma_i=\tau_i$ if $\fg\not= A_{2n}$). Let $\vec{\lambda}^\dagger$ be the $s$-tuple of dominant weights with $\lambda^\dagger_i\in D_{c, \sigma_i}$ associated to $\lambda_i\in D_{c, \gamma_i}$ via the bijection (\ref{wt_bij}). (The bijection $D_{c, \tau_i}\simeq D_{c, \sigma_i}$ between the diagram automorphism $\tau_i$ and the standard automorphism $\sigma_i$ for $\fg = A_{2n}$ is explicitly given by Lemma 
\ref{A_2n_lem}.)

For convenience, we fix a standard automorphism $\sigma$ corresponding to a generator of $\tilde{\Gamma}$.  

\begin{theorem}\label{dimensionformula}
With the same notation and Assumption  \ref{assump_of_thm}
as above for any finite group $\Gamma$, we further assume that the vanishing conjecture \ref{Teleman} holds. Then, the dimension of the twisted conformal blocks  
\[   N_{\bar{g},\Gamma}(\vec{\gamma};\vec{\lambda})=  \frac{  |T_c|^{\bar{g}-1+a }  }  { |T_c^\sigma|^{a} }    \sum_{t\in T_c^{\sigma, reg}/W^\sigma}  \frac{  \chi_{\vec{\lambda}^\dagger}(t)   \Delta_\sigma(t)^{a}  } {\Delta(t)^{\bar{g} -1+a }  }  ,\]
where $T$ is the maximal torus of $G$ with its Lie algebra the Cartan subalgebra $\fh$,  $\bar{g}$ is the genus of the quotient curve $\bar{\Sigma}=\Sigma/\Gamma$, and $\chi_{\vec{\lambda}^\dagger}(t) := \chi_{\lambda^\dagger_1}(t)\cdots \chi_{\lambda^\dagger_s}(t)$ with $\chi_{\lambda^\dagger_i}$ being the character of the irreducible representation of the group $G^{\sigma_i}$ with highest weight $\lambda^\dagger_i$. (For the notation $\Delta_\sigma$, see Theorem 
\ref{verlinde_Hong} and $\Delta= \Delta_{\text{identity}}$.)
\end{theorem}
\begin{proof} If $\tilde{\Gamma} = (1)$, the theorem follows from the reduction Corollary \ref{newcoro} and the classical Verlinde formula (cf.\,\cite[Theorem 4.2.19 and the identities (3) and (8) in its proof]{Kbook2}). So, we now assume that 
$\tilde{\Gamma} \neq  (1)$ in what follows.

 Set $A_c=T_c^{reg}/W $ and $A_c^\sigma=T_c^{\sigma, reg}/W^\sigma $. By Lemma \ref{lem_5.2}, $A^\sigma_c$ can be regarded as the subset of $\sigma$-invariants in $A_c$. 
 Clearly,  $A_c^\sigma= A_c^{\sigma_i}$ if $\langle  \sigma\rangle  = \langle  \sigma_i \rangle$, and $A_c^\sigma=A_c$ if $\sigma$ is trivial. 
 
 By Theorem \ref{red_thm}, we are reduced to the following situation: $\Gamma$ is cyclic  of order $2$ or $3$ (in particular, every non-trivial element is a generator of $\Gamma$) and the elements of $\Gamma$ act on $G$ by diagram automorphisms (possibly trivial), and the ramification type of $\Gamma$ action on $\Sigma$ is
  $\vec{\gamma}$ with each $\gamma_i$ acting on $G$ by a diagram automorphism $\tau_i$ and $\tau_i\not= 1$ exactly when $1\leq i \leq 2a$ ($a\geq 0$). The $s$-tuple $(\lambda_1, \lambda_2,\cdots, \lambda_{2a}, \mu_1, \cdots, \mu_{b})$ of dominant weights is attached to $\vec{p}$ with $\lambda_i\in D_{c, \gamma_i}$ and $\mu_j\in D_c$, 
  where $s=2a+b$. (Observe that $\mu_j := \lambda_{2a+j}$.) They satisfy $\tau_{2k-1}\tau_{2k}=1$ for any $1 \leq k\leq a$ as in Assumption \ref{assump_of_thm}.  Then, 
 \begin{align}
   N_{\bar{g},\Gamma}(\vec{\gamma};\vec{\lambda}, \vec{\mu}) &=   \sum_{\vec{\nu}\in D_c^a}   \left ( \prod_{k=1}^a N(\tau_{2k-1}; \lambda_{2k-1},\lambda_{2k}, \nu_k)  \right ) \cdot N_{\bar{g}}(\vec{\mu},\vec{\nu}^*  )   \nonumber \\
   &=   \sum_{\vec{\nu}}   \left ( \prod_{k=1}^a N(\sigma_{2k-1}; \lambda^\dagger_{2k-1},\lambda^\dagger_{2k}, \nu_k)  \right ) \cdot N_{\bar{g}}(\vec{\mu},\vec{\nu}^*  )   \nonumber \\
   &=\frac{  |T_c|^{\bar{g}-1 }  }  { |T_c^\sigma|^{a} }  \sum_{\vec{\nu}} \left ( \prod_{k=1}^a   \left(\sum_{t_k \in A^\sigma_c } \chi_{\lambda^\dagger_{2k-1}}(t_k) \chi_{\lambda^\dagger_{2k}}(t_k) \chi_{\nu_k }(t_k)\Delta_\sigma(t_k)\right) \right ) \cdot \left(\sum_{t\in A_c}   \chi_{\vec{\mu}}(t) \chi_{\vec{\nu}^*}(t) \Delta(t)^{1-\bar{g}}\right) \nonumber \\
   \label{eq_5.8_1}          &=\frac{  |T_c|^{\bar{g}-1 }  }  { |T_c^\sigma|^{a} }  \sum_{\vec{\nu}}  \sum_{\substack{t_1,\cdots, t_a \in A_c^\sigma\\  t \in  A_c }   }    \Delta(t)^{1-\bar{g}}  \chi_{\vec{\mu}}(t)    
    \prod_{k=1}^a  \left( \chi_{\lambda^\dagger_{2k-1}}(t_k) \chi_{\lambda^\dagger_{2k}}(t_k)    \chi_{\nu_k }(t_k) \chi_{\nu^*_k}(t)    \Delta_\sigma(t_k)\right),
 \end{align}
 where the first equality follows from Theorem \ref{g_red_thm_2}, the second equality follows from Theorem \ref{red_thm}, the third equality follows from Theorem \ref{verlinde_genus_0} and the usual Verlinde formula (cf.\,\cite[Theorem 4.2.19 and the identities (3) and (8) in its proof]{Kbook2}). 
 
Recall the following orthogonality relation (cf.\,\cite[Theorem 2 (2)]{Ho2}):
\[ \frac{1}{|T_c|}    \sum_{t\in  A_c  }  \chi_{\nu}(t) \chi_{\nu'^*}(t) \Delta(t)=\delta_{ \nu, \nu'} ,\]
where $\delta$ denotes the Kronecker symbol. Similarly, for any $t', t\in A_c$, the following orthogonality relation holds
\begin{equation}  \label{neweqn48}  \frac{1}{|T_c|}    \sum_{ \nu\in D_c  }  \chi_{\nu}(t') \chi_{\nu^*}(t) \Delta(t)=\delta_{t', t} . 
\end{equation}
Now, 
\begin{align}
 &\sum_{\vec{\nu}}  \sum_{\substack{t_1,\cdots, t_a \in A_c^\sigma\\  t \in  A_c }   }    \Delta(t)^{1-\bar{g}}  \chi_{\vec{\mu}}(t)    
    \prod_{k=1}^a  \left( \chi_{\lambda^\dagger_{2k-1}}(t_k) \chi_{\lambda^\dagger_{2k}}(t_k)    \chi_{\nu_k }(t_k) \chi_{\nu^*_k}(t)    \Delta_\sigma(t_k)\right)\nonumber\\
 & =  \sum_{\substack{t_1,\cdots, t_a \in A_c^\sigma\\  t \in  A_c }   } \,  \Delta(t)^{1-\bar{g}-a}  \chi_{\vec{\mu}} (t)
\left( \prod_{k=1}^a   \chi_{\lambda^\dagger_{2k-1}}(t_k) \chi_{\lambda^\dagger_{2k}}(t_k) \Delta_\sigma(t_k)\right) 
 \left( \sum_{\nu_1}\chi_{\nu_1^* }(t) \chi_{\nu_1}(t_1)  \Delta(t)\right) \dots  \left( \sum_{\nu_a}\chi_{\nu_a^* }(t) \chi_{\nu_a}(t_a)  \Delta(t)\right) \nonumber \\
 & = |T_c|^a\sum_{\substack{t_1,\cdots, t_a \in A_c^\sigma\\  t \in  A_c }   } \,  \Delta(t)^{1-\bar{g}-a}  \chi_{\vec{\mu}} (t)
 \prod_{k=1}^a  \left( \chi_{\lambda^\dagger_{2k-1}}(t_k) \chi_{\lambda^\dagger_{2k}}(t_k) \Delta_\sigma(t_k)\delta_{t, t_k}\right), \,\,\text{by using equation \eqref{neweqn48}} \nonumber \\
 \label{neweqn49}       & =\sum_{t\in A_c^\sigma}\,  \Delta(t)^{1-\bar{g}-a}  \chi_{\vec{\mu}} (t)  \chi_{\vec{\lambda}^\dagger} (t) (\Delta_\sigma(t))^a 
 |T_c|^a.
 \end{align}
 Combining the equations \eqref{eq_5.8_1}  and \eqref{neweqn49}, we get 
\[   N_{\bar{g},\Gamma}(\vec{\gamma};\vec{\lambda}, \vec{\mu})=     \frac{  |T_c|^{\bar{g}-1+a }  }  { |T_c^\sigma|^{a} }    \sum_{t\in A_c^\sigma}  \frac{  \chi_{\vec{\lambda}^\dagger}(t) \chi_{\vec{\mu}}(t)  \Delta_\sigma(t)^{a}  } {\Delta(t)^{\bar{g} -1+a }  }  . \]
   
This conclude the proof of the theorem. 
\end{proof}

\begin{remark}
By \cite[Prop 9.6]{DM},  the formulae for $S$-matrices described in \cite[$\S$ 5]{Ho2} and Lemma \ref{lem_5.2},  one  observes that the dimension formula in Theorem \ref{dimensionformula} agrees with the dimension formula in \cite[Theorem 1.2]{DM}.
\end{remark}

\end{document}